\documentclass[11pt,a4paper,english]{article}
\usepackage{jmlr2e}

\usepackage{amsmath}
\usepackage{amsfonts}
\usepackage{mathrsfs}
\usepackage{graphicx}
\usepackage{times}
\usepackage{bm}
\usepackage{natbib}

\usepackage[plain,noend]{algorithm2e}

\makeatletter
\renewcommand{\algocf@captiontext}[2]{#1\algocf@typo. \AlCapFnt{}#2} 
\def\@algocf@capt@plain{top}
\renewcommand{\algocf@makecaption}[2]{%
  \addtolength{\hsize}{\algomargin}%
  \sbox\@tempboxa{\algocf@captiontext{#1}{#2}}%
  \ifdim\wd\@tempboxa >\hsize
    \hskip .5\algomargin%
    \parbox[t]{\hsize}{\algocf@captiontext{#1}{#2}}
  \else%
    \global\@minipagefalse%
    \hbox to\hsize{\box\@tempboxa}
  \fi%
  \addtolength{\hsize}{-\algomargin}%
}
\makeatother

\newcommand{\bC}{{C}}
\newcommand{\bz}{{z}}

\newcommand{\bbeta}{{\beta}}
\newcommand{\bdelta}{{\delta}}

\newcommand{\bgamma}{{\gamma}}

\newcommand{\beps}{{\varepsilon}}
\newcommand{\bSigma}{{\Sigma}}

\newcommand{\bbetabar}{\mu}

\newcommand{\hbbetak}{{\widehat{\bbeta}_{k}}}
\newcommand{\hbbeta}{\widehat {\beta}}
\newcommand{\hbeta}{\widehat {\beta}}

\newcommand{\hbetabar}{\widehat {\mu}}
\newcommand{\hbbetabar}{\widehat {\bbetabar}}
\newcommand{\hbgamma}{\widehat {\gamma}}

\newcommand{\btheta}{{\theta}}

\newcommand{\bomega}{{\omega}}
\newcommand{\bOmega}{{\Omega}}
\newcommand{\btau}{{\tau}}

\newcommand{\bI}{{I}}
\newcommand{\by}{{y}}
\newcommand{\bY}{{y}}
\newcommand{\bX}{{X}}
\newcommand{\bx}{{x}}

\newcommand{\bcX}{{{\mathcal X}}}

\newcommand{\bcY}{{ {\mathscr{Y}}}}

\renewcommand{\P}{{\rm I}\kern-0.12em{\rm P}}
\newcommand{\1}{{\rm 1}\kern-0.28em{\rm I}}
\newcommand{\E}{{\rm I}\kern-0.12em{\rm E}}
\newcommand{\R}{{\rm I}\kern-0.14em{\rm R}}
\newcommand{\N}{{\rm I}\kern-0.14em{\rm N}}
\newcommand{\reals}{{\rm I}\kern-0.16em{\rm R}}

\newcommand{\p}{{\rm I}\kern-0.18em{\rm P}}

\DeclareMathOperator*{\argmin}{argmin}

\begin{document}

\title{Regression modeling on stratified data with the lasso}

\author{
\name{ E. Ollier$^{\star}$ and V. Viallon$^\ddag$} \\ \\
\addr{
$^\star$ Ecole Normale Sup\'erieure de Lyon 46 All\'ee d'Italie, 69364 Lyon Cedex 07, France \\
$^\ddag$Universit\'e de Lyon, F-69622, Lyon, France;
Universit\'e Lyon 1, UMRESTTE, F-69373 Lyon;
IFSTTAR, UMRESTTE, F-69675 Bron.\\
}
\email{vivian.viallon@univ-lyon1.fr}
\email{edouard.ollier@ens-lyon.fr}
}

\maketitle

\begin{abstract}
We consider the estimation of regression models on strata defined using a categorical covariate, in order to identify interactions between this categorical covariate and the other predictors. A basic approach requires the choice of a reference stratum.  We show that the performance of a penalized version of this approach depends on this arbitrary choice. We propose a refined approach that bypasses this arbitrary choice,  at almost no additional computational cost. Regarding model selection consistency, our proposal mimics the strategy based on an optimal and covariate-specific choice for the reference stratum. Results from an empirical study confirm that our proposal generally outperforms the basic approach in the identification and description of the interactions. An illustration is provided on gene expression data.  

\end{abstract}

\section{Introduction}\label{sec:Intro}

We consider the estimation of regression models when the population under study is stratified, as is standard in epidemiology and clinical research. For instance, when studying relapse after a primary breast cancer, it is now common to analyze various histological subtypes  at once \citep{voduc2010breast, rosner2013breast}. In order to accurately estimate the risk of relapse according to cancer subtype,  risk factors that interact with cancer subtype need to be identified, and the corresponding interactions need to be precisely described. In pharmacokinetics, it is often of interest to describe how parameters related to absorption and clearance depend on dosage and type of adjuvant. In  \cite{SAEM}, for example, non-linear mixed effect models are estimated on strata defined according to the treatment dosage and the type of adjuvant. In Section \ref{sec:RealData}, we study the association between the expression of the epidermal growth factor receptor gene, and those of 44 other transcription factors at eight time points of a differentiation process. We use a data set described in \cite{kouno2013temporal} where expression of these factors has been measured by profiling, for each time point, 120 different single cells. In this application, our aim is to describe how the association between the epidermal growth factor receptor gene and the other $44$ transcription factors varies over these eight strata.

In all these examples, the general objective is to study the relationship between a response variable $y\in\R$ and a vector of $p\geq 1$ predictors $\bx=(x_1, \ldots, x_p)\in\R^p$ over $K\geq 1$ strata defined according to a categorical covariate $Z$ of primary interest, such as cancer subtype. A key objective is  to determine how $Z$ modifies the effect of $\bx$ on $y$, that is to identify and describe interactions between predictors $\bx$ and the categorical effect modifier $Z$  \citep{GertheissTutz}. Identifying and taking account of these interactions is important even when the assessment of categorical effect modification is not the main objective. In particular, estimating one model per stratum or one model on all the strata pooled together generally leads to overfitting or underfitting, respectively. 

For simplicity, we mostly focus on linear regression models. Denote by $\bbeta^*_k=(\beta^*_{k,1}, \ldots, \beta^*_{k,p})\in\R^p$, for any $k\in[K]=\{1,\ldots, K\}$,  the parameter vector describing the association between $y$ and $\bx$ in the model corresponding to stratum $Z=k$. For any $j\in [p]$, denote by $d_j\in[K]$ the number of distinct values in the set $\{\beta^*_{1,j}, \ldots, \beta^*_{K,j}\}$. Further consider the partition $\{{\cal K}_j^{(1)}, \ldots, {\cal K}_j^{(d_j)}\}$ of $[K]$ such that $\beta^*_{k_1,j} = \beta^*_{k_2,j}$ if and only if  $(\beta^*_{k_1,j}, \beta^*_{k_2,j} )\in {\cal K}_j^{(d)}$ for some $d\in[d_j]$. The full description of how $Z$ modifies the effect of $x_j$ on $y$ relies on the identification of this partition. The total number of partitions of $[K]$ is the $K$th Bell number $B_K$, with, for instance, $B_5 = 52$. Standard statistical test procedures would require the comparison of  $(B_K)^p$ models for the identification of the $p$ partitions, so they are not well-suited to fully describe how $Z$ modifies the effect of $\bx$ on $y$.  

Penalized approaches have been advocated in this context, which can be seen as a special case of multi-task learning  \citep{EvgeniouPontil}. 
In particular, a version of the adaptive generalized fused lasso \citep{Tibshirani05, Viallon} has been shown to enjoy an asymptotic oracle property if $Kp$ does not grow with the sample size $n$  \citep{GertheissTutz, Tutz2}. If $Kp$ is fixed then this method identifies partitions $\{{\cal K}_j^{(1)}, \ldots, {\cal K}_j^{(d_j)}\}$ for all $j\in[p]$ with probability tending to one as $n\rightarrow\infty$, under mild assumptions. However, theoretical results in a non-asymptotic framework are still lacking for this method. Results obtained by \cite{Sharpnack} for generalized fused lasso estimates in the Gaussian means setting are not straightforward to extend to our context but they suggest that this method might only be able to identify the partitions under very particular settings. See Section \ref{sec:Ortho} for more details. 

An alternative strategy, which is standard in epidemiology, consists in first selecting a reference stratum, say the first one, then coding the other strata by indicators $\1(Z=k)$, and finally including them into the regression model, along with their interactions with the predictors. This corresponds to considering the decompositions $\bbeta^*_k = \bbeta^*_1 + \bdelta^*_k$, for all $k\in[K]$, with $\bdelta^*_1 = {0}_p$, the null vector in $\R^p$. Then, the lasso \citep{TibsLasso} can be used to identify null components in vectors $\bbeta^*_1$ and $\bdelta^*_k$, for $k\neq 1$. Of course, this strategy can generally identify only one element of each partition $\{{\cal K}_j^{(1)}, \ldots, {\cal K}_j^{(d_j)}\}$. However, a natural question is whether it is sparsistent, that is whether it can identify the sets $S^*_1 = \{j\in[p]: \beta^*_{1, j}\neq 0\}$ and $T^*_1 = \{(k, j)\in [K] \times [p]: \beta^*_{k,j}\neq \beta^*_{1,j}\}$ with high probability. If so, it could partially describe the effect of $Z$ on the association between $y$ and $\bx$, while results of \cite{Sharpnack} suggest that the method based on the fused lasso generally fails to fully describe it. 

In this article, we establish that the sparsistency of this basic approach depends on the choice of the reference stratum. We then present a refined approach which bypasses this arbitrary choice. Under linear regression models, we show that our proposal is sparsistent under conditions similar to those ensuring the sparsistency of the approach based on an optimal and covariate-specific choice of the reference stratum. In addition, our proposal can be implemented with available packages under a variety of models, at approximately the same computational cost as that of the basic approach. 


\section{Methods}\label{sec:Methods}

\subsection{Notations and setting } Methods are first presented in the linear regression model for ease of notation. Extensions to generalized linear models \citep{McC89} are briefly presented in Section \ref{sec:RewritingLasso}. 

For any positive integer $m\geq 1$, define $[m]=\{1,\ldots,m\}$. Let ${ 0}_m$ and ${1}_m$ be the vectors of size $m$ with components all equal to 0 and 1 respectively, and let ${I}_m$ be for the $(m\times m)$ identity matrix. For any vector $\bx=(x_1,\ldots,x_m)^T\in\R^m$, let ${\rm supp}(\bx)= \{j\in[m]: x_j \neq 0\}$ denote its support. We further set $\|\bx\|_q = (\sum_{j\in[m]} |x_j|^q)^{1/q}$ for any real number $q\in(0,\infty)$, $\|\bx\|_\infty = \max_j |x_j|$ and $\|\bx\|_0 = |{\rm supp}(\bx)|$, where $|E|$ is the cardinality of the set $E$. For any set $E\subseteq [m]$, let $\bx_E$ denote the vector of $\R^{|E|}$ with components $(x_j)_{j\in E}$. For any real matrix ${M}$, let $M_j$ be its $j$th column, and let ${M}_E$ be the sub-matrix made of columns $(M_j)_{j\in E}$. Further denote the smallest singular value of $M$ by $\Lambda_{\min}({M})$. Finally, $\1(\cdot)$ is the indicator function.

Denote the number of levels of variable $Z$, that is the number of strata, by $K\geq 1$. Let $n_k$ be the number of observations in stratum $k\in[K]$, and denote the total number of observations by $n=\sum_{k\in[K]}n_k$. For $k\in[K]$, further denote the response vector in stratum $k$ by $\by^{(k)} \in\R^{n_k}$. Similarly, let $\bX^{(k)}$ be the $({n_k}\times p)$ design matrix in stratum $k$. For all $k\in[K]$, we assume that $\by^{(k)} = \bX^{(k)} {\bbeta^*_k} + \beps^{(k)}$, with noise vector $\beps^{(k)} = (\varepsilon^{(k)}_1,\ldots,\varepsilon^{(k)}_{n_k})^T\in\R^{n_k}$. Vectors $\bbeta^*_k\in\R^p$ include the $Kp$ parameters to be estimated.

\subsection{Basic approach}\label{sec:RefLasso}

A basic approach consists in picking a reference stratum for any $j$, say $\ell_j$. Most often in practice, $\ell_j$ is chosen so that it does not depend on $j$; here we consider the most general version.  Introduce $\ell = (\ell_1, \ldots, \ell_p)\in[K]^p$, $\mu_\ell^* = (\beta^*_{\ell_1,1}, \ldots, \beta^*_{\ell_p, p})$ and $\bdelta^*_{k} = \bbeta^*_k -\mu_\ell^*$. The basic approach relies on the decomposition  $\bbeta^*_{k}  = \mu_\ell^*+ \bdelta^*_{k}$, for all $k\in[K]$, with $\delta^*_{\ell_j, j} = {0}$ \citep{GertheissTutz}. Following the lasso \citep{TibsLasso}, estimates of $\mu_\ell^*$ and $\bdelta^*_{k}$, $k\in[K]$, can be defined as
\begin{equation}
\hspace{-8pt}\argmin_{\substack{\bbetabar,\bdelta_1,\ldots,\bdelta_K\\ \delta_{\ell_j, j}=0 \ {\rm for\ all}\ j\in[p]}}\left\{ \sum_{k=1}^K \frac{\|\bY^{(k)} - \bX^{(k)}(\bbetabar+\bdelta_{k})\|_2^2 }{2n} + \lambda_1 \|\mu\|_1 + \sum_{k=1}^K \lambda_{2,k}\|\bdelta_{k}\|_1\right\},  \label{eq:RefLasso}
 \end{equation}
for appropriate non-negative $\lambda_1$ and $\lambda_{2,k}$. As will be seen in Sections \ref{sec:nonasymp} and \ref{sec:Ortho}, conditions ensuring the sparsistency of this approach  depend on the arbitrary choice of the vector of reference strata $\ell$. As a matter of fact, the underlying model dimension is $\|\mu_\ell^*\|_0 + \sum_{k\neq \ell} \|\bbeta^*_{k} - \mu_\ell^*\|_0$, which depends on $\ell$. The lowest  dimension, and then the best possible performance, is attained if $\mu_{\ell_j}^*$ is a mode of the collection of values $(0, \beta^*_{1,j}, \ldots, \beta^*_{K,j})$, ${\rm mode}(0, \beta^*_{1,j}, \ldots, \beta^*_{K,j})$, for all $j\in[p]$. For any $j\in[p]$,  such an optimal and covariate-specific reference stratum will be denoted by $\ell^*_j$ below. Because $\ell^*=(\ell^*_1, \ldots, \ell^*_p)$ is generally unknown, the corresponding optimal version of the basic approach cannot be implemented in practice.

\subsection{Our proposal}\label{Sec:ourproposal}
Our proposal aims at bypassing the arbitrary choice of the vector of reference strata $\ell$, while mimicking the optimal version of the basic approach. We consider an overparametrization involving $(K+1)p$ parameters, $\bbeta^*_{k} = \bbetabar^* + \bgamma^*_{k}$ for any $k\in[K]$. It generalizes the decomposition used in the basic approach in the sense that no coefficient is constrained to be zero here. For nonnegative values of $\lambda_1$ and $\lambda_{2,k}$'s, our proposal returns 
\begin{align}
\hspace{-8pt}(\hbbetabar,\hbgamma_1\ldots,\hbgamma_K)\in \argmin_{\bbetabar,\bgamma_1,\ldots,\bgamma_K}\left\{\sum_{k=1}^K \frac{\|\bY^{(k)} - \bX^{(k)}(\bbetabar+\bgamma_{k})\|_2^2 }{2n}
  + \lambda_1 \|\bbetabar\|_1 + \sum_{k=1}^K \lambda_{2,k}\|\bgamma_{k}\|_1\right\}.  \label{est_M1}
\end{align}
Working with a large  enough value for $\lambda_{2,r}$, for some $r\in[K]$, is equivalent to constraining $\hbgamma_r = {0}_p$, and then $\hbbetabar = \hbbeta_r$, and reduces to the basic approach with $\ell=(r, \ldots, r)$. More generally, setting $\btau = (\tau_1,\ldots,\tau_K)$ with $\tau_k = \lambda_{2,k}/\lambda_1$, and defining the shrunk and $\btau$-weighted version of the median of $(b_1,\ldots, b_K)$ as ${\rm WSmedian}(b_1,\ldots, b_K; \btau) = \argmin_b (|b| + \sum_{k\in[K]} \tau_k |b_k - b|)$, it is easy to see that $\hbetabar_j \in {\rm WSmedian}(\hbeta_{1,j},\ldots, \hbeta_{K,j}; \btau)$. In other words, for any particular value of the $\lambda_{2,k}/\lambda_1$ ratios,  our approach encourages solutions $(\hbbeta_{1} = \hbbetabar + \hbgamma_1,\ldots, \hbbeta_{K}= \hbbetabar + \hbgamma_K)$ with a sparse vector $\hbbetabar$ and sparse vectors of differences $\hbgamma_k = \hbbetak - \hbbetabar$, and with the overall effect of the $j$th covariate $\hbetabar_j$ defined as ${\rm WSmedian}(\hbeta_{1,j},\ldots, \hbeta_{K,j}; \btau)$. 

Moreover, working with a large enough value for $\lambda_1$  is equivalent to constraining $\hbbetabar = {0}_p$, and our approach then reduces to $K$ independent lasso's. In contrast, working with large enough $\lambda_{2,k}$ values is equivalent to constraining $\hbbeta_k = \hbbetabar$  for all $k$, and our approach then reduces to one lasso run on all the strata pooled together. 

\subsection{Rewriting as a lasso on a transformation of the original data}\label{sec:RewritingLasso}

Our proposal reduces to the lasso on a simple transformation of the original data, just as the basic approach does. Set $\bcY = (\bY^{(1)^T}, \ldots, \bY^{(K)^T})^T$ the vector containing the $n$ observations of the response variable. For any $k\in[K]$, introduce $P^{(k)}_\ell = \{j\in[p]: k\neq \ell_j\}$, with $\ell_j$ still denoting the reference stratum chosen for covariate $j$ in the basic approach. Note that $\sum_k |P^{(k)}_\ell|=(K-1)p$ and set $\tilde{\bX}_\ell^{(k)}=\bX^{(k)}_{P^{(k)}_\ell}$. Now introduce
$$ 
\bcX_{\ell} =  \left( 
\begin{array}{c c c c}
\bX^{(1)} &\tilde{\bX}_\ell^{(1)}/\tau_1 & \hdots &{0} \\
\vdots & \vdots & \ddots & \vdots \\
\bX^{(K)} &{0}&\hdots  & \tilde{\bX}_\ell^{(K)}/\tau_K
\end{array}
\right), \quad \quad\bcX_0 =\left( 
\begin{array}{c c c c }
\bX^{(1)} &\bX^{(1)}/\tau_1& \hdots &{0}  \\
\vdots & \vdots & \ddots & \vdots \\
\bX^{(K)} &{0}&\hdots  &  \bX^{(K)}/\tau_K 
\end{array}
\right),
$$
Criteria to be minimized in (\ref{eq:RefLasso}) and (\ref{est_M1}) reduce to the lasso 
\begin{equation}
 \frac{1}{2n} \|\bcY - \bcX\btheta\|_2^2 + \lambda_1\|{\btheta}\|_{1}, \label{eq:M1_Lasso}
\end{equation}
with $\bcX$ set to  either $\bcX_\ell$, for the basic approach, or $\bcX_0$, for our proposal, and $\btheta$ a vector of $\R^{Kp}$ or $\R^{(K+1)p}$ as appropriate. Therefore, our proposal comes at almost no additional computational cost compared to the basic approach.

This rewriting as a lasso extends to generalized linear models, Cox models, etc., and makes our proposal directly implementable using the glmnet  package of \cite{glmnet}, for instance. Considering logistic models, set ${\mathcal L}_{\rm logistic}(\by,\bz) =  \sum_{i\in[n]} y_iz_i - \log(1+e^{z_i})$ for any $\by\in \{0,1\}^n$ and $\bz\in\R^n$. The criterion of our proposal writes as the logistic lasso,$ -{\mathcal L}_{\rm logistic}(\bcY,\bcX_0\btheta) + \lambda_1\|\btheta\|_{1}.$ In addition, the glmnet package can benefit from the sparse structure of $\bcX_0$.  

\subsection{Sparsistency}\label{sec:nonasymp}

From now on we assume that $\bar\beta^*_j = {\rm mode}(0, \beta^*_{1,j}, \ldots, \beta^*_{K,j})$ is uniquely defined for all $j\in[p]$, for ease of notation. For any vector of reference strata $\ell=(\ell_1, \ldots, \ell_p)\in[K]^p$, introduce the sets $S_{\ell}=\{j\in[p]: \beta^*_{\ell_j, j}\neq 0\}$ and $T_{\ell}=\{(k,j)\in[K]\times[p]: \beta^*_{k,j}\neq \beta^*_{\ell_j, j}\}$. Further define  $\btheta_\ell^{*}=({\mu_\ell^{*T}}, \tau_2{\bgamma_{\ell, 2}^{*T}}, \ldots,\tau_K\bgamma_{\ell,K}^{*T})^T\in\R^{Kp}$, with $\mu_{\ell,j}^{*}=\beta^*_{\ell, j}$ and $\bgamma_{\ell, k}^{*} = (\bbeta^*_k - \mu_\ell^{*})_{P^{(\ell)}_k}$. The basic approach is sparsistent if it identifies these two sets, or, equivalently, the set $J_{\ell} = {\rm supp}(\theta_\ell^{*})$, with high probability. Now, consider an optimal vector of reference strata $\ell^*=(\ell^*_1, \ldots, \ell^*_p)$ such that $\beta^*_{\ell^*_j, j}=\bar\beta^*_j$ for all $j\in[p]$.  Define $\btheta_0^{*}=(\bbetabar_{\ell^*}^{*T}, \tau_1\bgamma_{0,1}^{*T}, \ldots, \tau_K\bgamma_{0,K}^{*T})^T\in\R^{(K+1)p}$ where $\bgamma^{*}_{0,k} = \bbeta^*_k - \bbetabar_{\ell^*}^{*}$. We will say our proposal is sparsistent if it identifies $S_{\ell^*}$ and $T_{\ell^*}$,  or, equivalently, the set $J_0={\rm supp}(\theta_0^{*})$, with high probability. 

For the lasso to be sparsistent, a sufficient and almost necessary condition on the design matrix is the irrepresentability condition \citep{ZhaoYu, Wainwright2009}. Consider the general formulation (\ref{eq:M1_Lasso}) of the lasso and denote by $\btheta^*$ the true value of the parameter vector to be estimated. Defining $J^*={\rm supp}(\btheta^*)$, the matrix $\bcX$ fulfills the irrepresentability condition, with respect to $J^*$, if $\Lambda_{\min}(\bcX_{J^*}^T\bcX_{J^*})\geq C_{\min}$ for some $C_{\min}>0$ and $\max_{j\notin J^*}\|(\bcX_{J^*}^T\bcX_{J^*})^{-1}\bcX^{T}_{J^*}{\cal X}_j \|_1<1$. 
Using the vector of reference strata $\ell$ in the basic approach, the irrepresentability condition of matrix $\bcX_\ell$  writes as $(IC)_\ell$, while in our approach, the irrepresentability condition of matrix $\bcX_0$  writes as $(IC)_0$:
\begin{eqnarray*}
&(IC)_\ell& \quad \Lambda_{\min}(\bcX_{\ell J_\ell}^T\bcX_{\ell J_\ell})\geq C_{\ell}>0 \  {\rm and}\ c_\ell = \max_{j\notin J_\ell}\|(\bcX_{\ell J_\ell}^T\bcX_{\ell J_\ell})^{-1}\bcX^{T}_{\ell J_\ell}{\cal X}_{\ell j} \|_1<1, \\
&(IC)_0&\quad \Lambda_{\min}(\bcX_{0 J_0}^T\bcX_{0 J_0})\geq C_{0}>0 \ {\rm and}\ c_0 = \max_{j\notin J_0}\|(\bcX_{0 J_0}^T\bcX_{0 J_0})^{-1}\bcX^{T}_{0 J_0}{\cal X}_{0 j} \|_1<1.
\end{eqnarray*}
The comparison of $(IC)_{\ell^*}$ and $(IC)_0$ is of particular interest. First, because $\theta^*_{\ell^* J_{\ell^*}} = \theta^*_{0 J_0}$, we have ${\cal X}_{\ell^* J_{\ell^*}} =  \bcX_{0 J_0}$ and $\Lambda_{\min}(\bcX_{\ell^* J_{\ell^*}}^T\bcX_{\ell^* J_{\ell^*}}) = \Lambda_{\min}(\bcX_{0 J_0}^T\bcX_{0 J_0})$. Second, the maxima in the definitions of $c_0$ and $c_{\ell^*}$ are taken over $J_0^c$ and $J_{\ell^*}^c$, respectively, with $|J_0^c| = p + |J_{\ell^*}^c|$. Indeed, columns corresponding to $\gamma^*_{0,\ell^*_j,j}$, for $j\in[p]$, are present in matrix $\bcX_0$ while they are absent from $\bcX_{\ell^*}$. Moreover, the corresponding indexes belong to $J^c_0$ since $\gamma^*_{0, \ell^*_j,j}=0$.  Therefore, the only difference between $(IC)_0$ and $(IC)_{\ell^*}$ comes from the fact that $c_0\geq c_{\ell^*}$, and $(IC)_0$ is only slightly stronger than $(IC)_{\ell^*}$. Focusing on the case of balanced strata and orthogonal designs in each stratum, Lemma \ref{Lemma_IC_Ortho} in Section \ref{sec:Ortho} states that these two conditions are identical in this particular case. It further explicitly relates them to the ratios $\tau_k = \lambda_{2,k}/\lambda_1$ and to the maximum level of heterogeneity that is allowed among the collections of values $(\beta^*_{1,j}\ldots, \beta^*_{K,j})$, for all $j\in[p]$. In addition, they are shown to be generally weaker than $(IC)_\ell$ for the choice $\ell=(r, \ldots, r)$, for some $r\in[K]$.

We can now state Theorem \ref{theo:General}, according to which our proposal identifies $S_{\ell^*}$ and $T_{\ell^*}$ under nearly the same assumptions as those required by the optimal version of the basic approach if the $\ell^*_j$'s were given in advance. More precisely, besides the fact that $(IC)_0$ is generally a little stronger that $(IC)_{\ell^*}$, the only difference lies in the terms $(\lambda^{(1)}_1,\beta^{(1)}_{\min})$ and $(\lambda^{(0)}_1, \beta^{(0)}_{\min})$, where $K+1$ replaces $K$. Our result is a consequence of Theorem 1 in \cite{Wainwright2009}. 

\begin{theorem}\label{theo:General}
Assume that the noise variables $(\varepsilon^{(k)}_i)_{i\in[n_k], k\in [K]}$ are independent and identically distributed centered sub-Gaussian variables with parameter $\sigma>0$. Further assume that $n_k^{-1/2}\|X_j^{(k)}\|_2\leq 1$ for all $(k,j)\in[K]\times[p]$.  Introduce $\tau_0>0$ and set $\tau_k=\tau_0(n_k/n)^{1/2}$ for all $k\in[K]$. If $(IC)_{\ell^*}$ holds then set $\gamma_1 = 1 - c_{\ell^*}$. Further set $\gamma_0 = 1 - c_0$ if $(IC)_0$ holds. 
For $\eta\in\{0,1\}$, introduce 
$$ \lambda^{(\eta)}_1 > \frac{2}{\gamma_\eta\min(1,\tau_0)}\left\{\frac{2\sigma^2 \log((K+\eta)p))}{n}\right\}^{1/2}\!\!\!\!\!,\quad  \beta^{(\eta)}_{\min}= \lambda^{(\eta)}_1\left\{\frac{(|S_{\ell^*}| + |T_{\ell^*}|)^{1/2}}{C_{\ell^*}} + 4\frac{\sigma} {C_{\ell^*}^{1/2}}\right\}.$$

If $(IC)_{\ell^*}$ holds then solutions $\widehat{\btheta}_{\ell^*}$ of (\ref{eq:M1_Lasso}) with $\bcX=\bcX_{\ell^*}$ and $\lambda_1=\lambda_1^{(0)}$ as above are such that, with probability at least $1 - 4\exp(-a_0n\lambda_1)$ for some constant $a_0>0$:
$(i)$ $\widehat{\btheta}_{\ell^*}$ is uniquely defined, $(ii)$ $\widehat{\btheta}_{\ell^* J_{\ell^*}^c} = {0}_{|{J_{\ell^*}^c} |}$, and $(iii)$ $\|\widehat{\btheta}_{\ell^*J_{\ell^*}} -\btheta^*_{\ell^* J_{\ell^*}} \|_\infty\leq \beta^{(0)}_{\min}$. If, in addition,  $|\bar\beta^*_{j}|> \beta^{(0)}_{\min}$ for all $ j\in S_{\ell^*}$ and $|\beta^*_{k,j} - \bar\beta^*_{j}| > \beta^{(0)}_{\min}/\tau_k$ for all $(k,j)\in T_{\ell^*}$, then $J_{\ell^*}$, hence both $S_{\ell^*}$ and $T_{\ell^*}$, are perfectly identified with probability at least $1 - 4\exp(-a_0n\lambda_1^2)$.

If $(IC)_{0}$ holds then solutions $\widehat{\btheta}_{0}$ of (\ref{eq:M1_Lasso}) with $\bcX=\bcX_{0}$ and $\lambda_1=\lambda_1^{(1)}$ as above are such that, with probability at least $1 - 4\exp(-a_1n\lambda_1)$ for some constant $a_1>0$:
$(i)$ $\widehat{\btheta}_{0}$ is uniquely defined, $(ii)$ $\widehat{\btheta}_{0 J_0^c} = {0}_{|{J_0^c} |}$, and $(iii)$ $\|\widehat{\btheta}_{0 J_0} -\btheta^*_{\ell^* J_{\ell^*}} \|_\infty\leq \beta^{(1)}_{\min}$. If, in addition,  $|\bar\beta^*_{j}|> \beta^{(1)}_{\min}$ for all $ j\in S_{\ell^*}$ and $|\beta^*_{k,j} - \bar\beta^*_{j}| > \beta^{(1)}_{\min}/\tau_k$ for all $(k,j)\in T_{\ell^*}$, then $J_0$,  hence both $S_{\ell^*}$ and $T_{\ell^*}$, are perfectly identified  with probability at least $1 - 4\exp(-a_1n\lambda_1^2)$.
\end{theorem}

This result especially confirms that it is harder to identify $T_{\ell^*}$ than $S_{\ell^*}$, in the sense that heterogeneities have to be at least $|\beta^*_{k,j} - \bar\beta^*_j|> (n/n_k)^{1/2}\beta_{\min}^{(\eta)}/\tau_0$ for $(k,j)\in T_{\ell^*}$, while $|\bar\beta^*_j|$ has only to be greater than $\beta^{(\eta)}_{\min}$ for $j\in S_{\ell^*}$.

\subsection{The orthogonal and balanced case}\label{sec:Ortho}
It is instructive to inspect in more detail the simple setting where $n_k=n/K$ and $(\bX^{(k)^T}\bX^{(k)})/n_k = \bI_{n_k}$ for all $k\in[K]$. This orthogonality assumption does not make matrices $\bcX_{\ell}$ and $\bcX_0$ orthogonal and is therefore not sufficient for the irrepresentability condition. 

For any vector of reference strata $\ell$ and all $j\in[p]$, define $K^*_{\ell, j}=\{k\in[K]:\beta^*_{k,j}= \beta^*_{\ell_j,j}\}$. Set ${\cal D}_{\ell,0} = \max_{j\notin S_\ell } |\{k\in[K]:\beta^*_{k,j}\neq \beta^*_{\ell_j,j}\}|$ if $S_\ell \neq [K]$ and $0$ otherwise and ${\cal D}_{\ell,1}=\max_{j\in S_\ell} |\{k\in[K]:\beta^*_{k,j}\neq \beta^*_{\ell_j,j}\}|$ if $S_\ell\neq \emptyset$ and $-\infty$ otherwise. For all $k\in[K]$, we set $\tau_k = \tau_0 K^{-1/2}$,  for some $\tau_0>0$.

\begin{lemma}\label{Lemma_IC_Ortho}

The matrix $\bcX_\ell$ fulfills the irrepresentability condition if and only if 
\begin{equation*}
(sIC)_\ell\quad\quad0\leq \frac{K^{1/2}}{K - 2 {\cal D}_{\ell,1}}< \tau_0<\frac{K^{1/2}}{{\cal D}_{\ell,0}}\cdot
\end{equation*}
The matrix $\bcX_0$ fulfills the irrepresentability condition if and only if 
\begin{equation*}
(sIC)_0\quad\quad0\leq \frac{K^{1/2}}{K - 2 {\cal D}_{\ell^*,1}}< \tau_0<\frac{K^{1/2}}{{\cal D}_{\ell^*,0}}\cdot
\end{equation*} 

\end{lemma}

Conditions $(sIC)_0$ and $(sIC)_{\ell^*}$ in Lemma \ref{Lemma_IC_Ortho} are identical. In the orthogonal and balanced case, the two sets of assumptions required by our proposal and the optimal version of the basic approach to identify the sets $S_{\ell^*}$ and $T_{\ell^*}$ are therefore identical, except for the terms $(\lambda^{(1)}_1,\beta^{(1)}_{\min})$ and $(\lambda^{(0)}_1, \beta^{(0)}_{\min})$ where $K+1$ replaces $K$, as in Theorem \ref{theo:General} above. 

In addition,  $(sIC)_{\ell^*}$, or equivalently $(sIC)_{0}$, imposes that $2 {\cal D}_{\ell^*,1} + {\cal D}_{\ell^*,0} < K$. In particular, if  ${\cal D}_{\ell^*,1}  = {\cal D}_{\ell^*,0} = {\cal D}_{\ell^*}$, this implies that ${\cal D}_{\ell^*}<K/3$. The irrepresentability conditions $(sIC)_{\ell^*}$ and $(sIC)_{0}$ are therefore directly related to the maximum level of heterogeneity among the values $(\beta^*_{1,j}, \ldots, \beta^*_{K,j})$.  Similarly, $(IC)_\ell$ imposes $2 {\cal D}_{\ell,1} + {\cal D}_{\ell,0} < K$, which is generally a stronger constraint. For simplicity, consider the situation where $\ell=(r, \ldots, r)$ for some $r\in[K]$, which is a common choice in practice. Without loss of generality, set $r=1$. Then $(IC)_\ell$ entails that, for each $j\in[p]$, the effect of the $j$th covariate on most strata is $\beta^*_{1,j}$, while $(sIC)_{\ell^*}$ and $(sIC)_{0}$ only entail that, for each $j\in[p]$, the effect of the $j$th covariate on most strata is ${\rm mode}(0, \beta^*_{1,j}, \ldots, \beta^*_{K,j})$.  Moreover, if $\ell=(1,\ldots, 1)$ and $\{1\}\notin \cap_{j\in[p]}K^*_{\ell^*,j}$ then we have $T_\ell \neq T_{\ell^*}$ and, possibly, $S_\ell \neq S_{\ell^*}$: the identification of $T_\ell$ and $S_\ell$ is both less interesting and less likely and our proposal should be preferred over the basic approach. 



To recap, our non-asymptotic analysis shows that the partial description of the categorical effect modification due to $Z$ through decompositions of the type $\bbeta^*_k = \bbetabar^* + \bgamma^*_k$ is guaranteed only when the level of heterogeneity is not too high, that is when the number of nonzero components in vectors $(\gamma^*_{1,j}, \ldots,\gamma^*_{K,j})$, for $j\in[p]$, is not too high. In particular, the lowest level of heterogeneity is attained for the choice $\bbetabar^*_j = {\rm mode}(0, \beta^*_{1,j}, \ldots, \beta^*_{K,j})$. Our proposal is able to target this optimal decomposition and is sparsistent under nearly the same assumptions as those required for the optimal version of the basic approach.

\subsection{Connection with the generalized fused lasso}\label{sec:Fused}

The fact that the level of heterogeneity must be low to ensure the sparsistency of our proposal as well as the sparsistency of the optimal version of the basic approach has connections with other results in the literature. Consider for instance the approach mentioned in Section \ref{sec:Intro}, which was proposed by \cite{GertheissTutz}; see also \cite{Tutz2} and \cite{Viallon}. This is based on a fusion penalty that encourages similarities among solutions $\hbbeta_k\in\R^p$, $k\in[K]$. More precisely, for appropriate non-negative $\lambda_1$ and $\lambda_{2}$'s, it returns estimators $\hbbeta_k$ defined as minimizers of the criterion
\begin{equation}
 \sum_{k\in[K]} \frac{\|\bY^{(k)} - \bX^{(k)}\bbeta_{k}\|_2^2}{2n}  + \lambda_1 \sum_{k\in[K]} \|\bbeta_k\|_1 + \lambda_{2}\sum_{k_1< k_2} \|\bbeta_{k_1} - \bbeta_{k_2}\|_1.   \label{eq:CliqueFused}
 \end{equation}
 This criterion is that of the generalized fused lasso, where the graph used in the penalty is made of $p$ cliques of size $K$ \citep{Viallon}. The $j$th clique corresponds to the $j$th predictor and connects all the components in $(\beta_{1,j}, \ldots, \beta_{K,j})$: the $K(K-1)/2$ differences $|\beta_{k_1, j}-\beta_{k_2, j}|$, for $k_1<k_2$, appear in the penalty term. 
 
Both the basic approach and our proposal are related to generalized fused lasso estimates too. In particular, consider criterion (\ref{eq:RefLasso}) with the optimal choice $\ell^*$ for the vector of reference strata. It can be seen as a version of a generalized fused lasso, with a graph made of $p$ star-graphs of size $K$ instead of $p$ cliques: for each $j\in[p]$, only the $K-1$ differences $|\beta_{k,j}-\beta_{\ell^*_j,j}|$ appear in the penalty term, for $k\neq \ell^*_j$ and $\ell^*_j\in[K]$ fixed. 

Non-asymptotic analyses of the sparsistency of generalized fused lasso estimates are scarce in the literature. \cite{Sharpnack} study them in the normal means setting, which can be seen as a special case of the stratified linear regression considered here. \cite{Sharpnack} establish that generalized fused lasso estimates  are sparsistent only if the graph used in the fused penalty is in good agreement with the true structure of the vector of parameters; see also \cite{QianJia}. Although it is not straightforward to extend to our case, these results suggest that estimates derived from (\ref{eq:CliqueFused}) can only be sparsistent if the level of heterogeneity is not too high. 

Our results precisely quantify the maximum level of heterogeneity above which a version of generalized fused lasso estimates, based on star-graphs, can attain sparsistency in stratified regression, in the balanced and orthogonal case. Because star-graphs are less connected than cliques, it is likely that sparsistency for clique-based estimates, such as those minimizing criterion  (\ref{eq:CliqueFused}), requires an even lower maximum level of heterogeneity. That being said, sparsistency for clique-based estimates refers to the full identification of the partitions $\{{\cal K}_j^{(1)}, \ldots, {\cal K}_j^{(d_j)}\}$, for $j\in[p]$. For the basic approach, its optimal version and our proposal, sparsistency refers to the identification of one element of this partition only, say ${\cal K}_j^{(1)}$, and its complementary $[K]\setminus {\cal K}_j^{(1)}$. 

In the asymptotic regime, assuming that $Kp$ is fixed and $n_k/n \rightarrow \rho_k$ for some $\rho_k\in(0,1)$, for all $k\in[K]$, oracle properties have been derived for adaptive versions of clique-based estimates under mild assumptions \citep{GertheissTutz}: in particular, no assumption regarding the level of heterogeneity is required to ensure perfect recovery of the full partition $\{{\cal K}_j^{(1)}, \ldots, {\cal K}_j^{(d_j)}\}$, for all $j\in[p]$. Similar results are easily derived for adaptive versions of the basic approach for instance. In view of (\ref{eq:M1_Lasso}),  we can apply Theorem 2 of \cite{Zou06} to show that an adaptive version of the basic approach enjoys an oracle property too, without having to assume any irrepresentability condition. Here again, no assumption regarding the maximum level of heterogeneity is required, but the identification of only one element of the partition  is guaranteed for all $j$. 

To recap, clique-based estimates are optimal and should be preferred over our proposal or the basic approach in the asymptotic regime, assuming that $Kp$ is fixed and $n_k/n \rightarrow \rho_k$ for some $\rho_k\in(0,1)$, for all $k\in[K]$. In a non-asymptotic setting, results for clique-based estimates are still lacking, while conditions ensuring the sparsistency of our proposal and the basic approach are established in the present article. In the following simulation study, we especially compare the clique-based strategy and our proposal on finite samples.

\section{Simulation Study}\label{sec:Simulations}

Theorem \ref{theo:General} states that our proposal and the optimal version of the basic approach perform similarly with regard to the identification of $S_{\ell^*}$ and $T_{\ell^*}$ for appropriate values of $\lambda_1$ and $\tau_0$, under technical assumptions on the design matrices. The main objective of this simulation study is to assess the empirical performance of our proposal under general designs, and for $\lambda_1$ and $\tau_0$ selected by 5-fold cross-validation. Comparisons are made with the basic approach with the choice $\ell=(1, \ldots, 1)$ as well as an optimal choice $\ell^*$. Clique-based estimates are also considered. 

\begin{figure}[h]
\begin{center}
\includegraphics[scale=0.38]{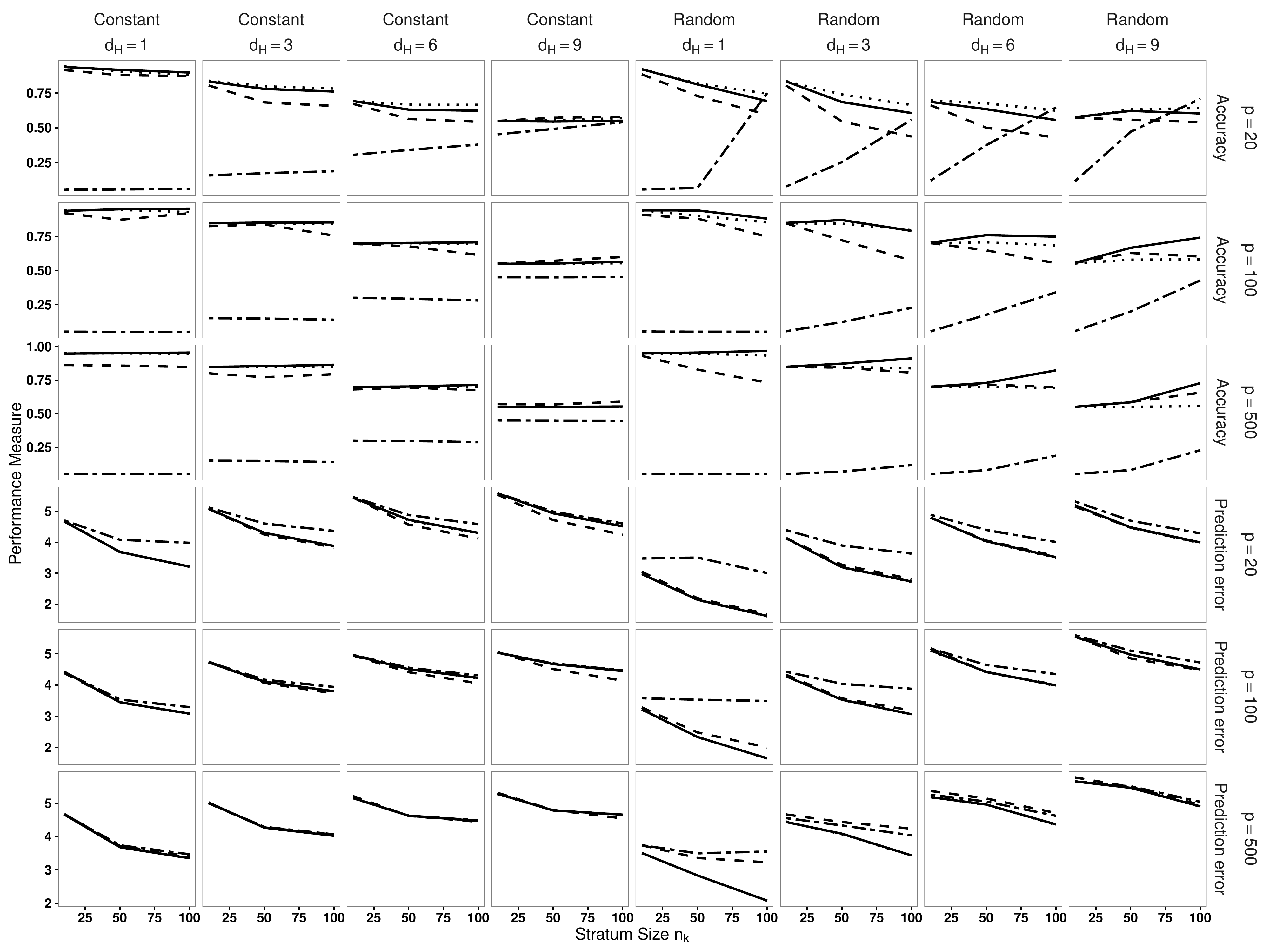}
\end{center}\caption{ Results from the simulation study. The top three panels show the accuracy regarding the recovery of the set $T^{*}_{1, P_0}=\{(k,j)\in[K]\times P_0: \beta^*_{k,j}\neq \beta^*_{1,j}\}$ for the basic approach and the set $T^*_{P_0}=\{(k,j)\in[K]\times P_0: \beta^*_{k,j}\neq \beta^*_{\ell^*_j,j}\}$ for the other approaches. The higher, the better. The bottom three panels illustrate the prediction error; the lower, the better. Results are presented for both the constant and random $\delta^*_{k,j}$ cases. All results correspond to averages over 50 replicates in each configuration. 
Solid line: our proposal. Dotted line: optimal version of the basic approach. Dash-dot line: basic approach. Dashed line: clique-based approach.}\label{fig:Accs}
\end{figure}

We set $K=20$ and take $n_k\in\{10, 50, 100\}$ and $p\in\{20, 100, 500\}$. For each $k\in[K]$, rows of the design matrix $\bX^{(k)}$ are drawn from an ${\cal N}({0}_p, \bSigma)$ distribution, with $\bSigma$ the $(p\times p)$ Toeplitz matrix with element $(i,j)$ equal to $0.5^{|i-j|}$. We then randomly select a subset $P_0\subset [p]$ of size $20$ and set  $\beta^*_{k,j} = 0$ for all $j\notin P_0$ and $k\in[K]$. As for the values $(\beta^*_{1,j}, \ldots, \beta^*_{K,j})$ for $j\in P_0$, we consider four levels of heterogeneity $d_H$. More precisely, for any given $d_H\in\{1, 3, 6, 9\}$, we set $\beta^*_{k,j} = 1$ for $k> d_H$, and  $\beta^*_{k,j} = 1 + \delta^*_{k,j}$ for $k\leq d_H$ for $10$ indexes $j$ randomly selected in $P_0$. For the other $10$ indexes in $P_0$, we set $\beta^*_{k,j} = 1$ for $k\leq d_H$, and  $\beta^*_{k,j} = 1 + \delta^*_{k,j}$ for $k> d_H$. We further consider two cases for the $\delta^*_{k,j}$ values:  they are either constantly set to $K^{1/2}$ or drawn from the uniform distribution on $[K^{1/2}/2, 2K^{1/2}]$ and then multiplied by $\pm 1$ (with probability $1/2$). When the $\delta^*_{k,j}$'s are constant, the collection of values $(\beta^*_{1,j}, \ldots, \beta^*_{K,j})$ for each covariate $j\in P_0$ is made of  two groups of distinct values, of sizes $K-d_H$ and $d_H$. This situation should favor clique-based estimates. When $\delta^*_{k,j}$ is random, the collection of values $(\beta^*_{1,j}, \ldots, \beta^*_{K,j})$ for each covariate $j\in P_0$ is made of $d_H+1$ groups, one of size $K-d_H$, and the other $d_H$ of size $1$. This situation should favor our proposal and the optimal version of the basic approach. Observations of the response variable are then generated according to $\by^{(k)} = \bX^{(k)} \bbeta^*_k + \beps^{(k)}$, with each component of $\beps^{(k)}$ drawn from an ${\cal N}(0,\sigma^2)$ distribution. The variance $\sigma^2$ is set to $\sum_{k\in[K]} \|\bX^{(k)}\bbeta^*_k\|_2^2/n$, giving overall signal-to-noise ratio equal to 1. For each particular choice of $n_k\in\{10, 50, 100\}$, $p\in\{20, 100, 500\}$ and $d_H\in\{1, 3, 6, 9\}$, and for both the random and constant choice for the $\delta^*_{k,j}$'s, we generate 50 replicates of data $(\bX^{(k)}, \by^{(k)})$, $k\in[K]$. Our results correspond to averages over these 50 replicates; see Figure \ref{fig:Accs}.

In all the configurations, $\beta^*_{20, j} = {\rm mode}(0, \beta^*_1, \ldots, \beta^*_{20,j})$ for all $j\in[p]$. We then set $\ell^*_j=20$ for all $j\in[p]$ for the optimal version of the basic approach. On the other hand, $\beta^*_{1, j} \neq {\rm mode}(0, \beta^*_1, \ldots, \beta^*_{20,j})$ for all $j\in P_0$. Under this setting, which is of course extreme, the comparison between the results obtained using either $\ell$ or $\ell^*$ for the reference strata allows a precise description of the impact of the reference stratum on the performance of the basic approach. Top panels of Figure \ref{fig:Accs} present results regarding the identification of the sets $T^*_{P_0}=\{(k,j)\in[K]\times P_0: \beta^*_{k,j}\neq \beta^*_{\ell^*_j,j}\}$ for the optimal version of the basic approach, our proposal and clique-based estimates, and that of the set $T^{*}_{1, P_0}=\{(k,j)\in[K]\times P_0: \beta^*_{k,j}\neq \beta^*_{1,j}\}$ for the basic approach. Here, we only consider covariates in $P_0$ because they are those that are the most differently accounted for by the four approaches we compare. 

In the constant $\delta^*_{k,j}$ case, our empirical results clearly illustrate our theoretical ones. First, our proposal performs nearly as well as the optimal version of the basic approach. Second, the lower $d_H$, the better they perform, as expected since $d_H = {\cal D}_1$ and ${\cal D}_0 = 0$ here. Third, it is more difficult to recover $T^{*}_{1, P_0}$ than $T^{*}_{P_0}$, which is also expected since $ {\cal D}^{(1)}_1 = K-d_H > {\cal D}_1$. For the same reason however, as $d_H$ increases,  $T^{*}_{1, P_0}$ is easier to recover. A nice symmetry appears between the performance of our proposal and the optimal version of the basic approach on the one hand and the basic approach using $\ell=(1, \ldots, 1)$, on the other hand. In addition, the recovery of the sets $T^*_{P_0}$ and $T^{*}_{1, P_0}$ is only marginally affected by $p$. 
Results in the random $\delta^*_{k,j}$ case are mostly consistent with those in the constant case, except for the recovery of  $T^{*}_{1, P_0}$ which is mostly due to the fact that $T^{*}_{1, P_0} = [K-1]\times P_0$ irrespective of $d_H$ in this random case. Finally, the clique-based strategy performs similarly to, or a little worse than, our proposal with respect to this criterion. 

The bottom panels of Figure \ref{fig:Accs} present results regarding  $\log(\sum_{k\in [K]} \|\bX^{(k)}(\bbeta^*_k - \hbeta_k)\|_2^2/n)$  \citep{Dalalyan14}, which is a measure of the prediction error. Overall, our proposal and the optimal version of the basic approach perform similarly. They both outperform the basic approach, especially in the random $\delta^*_{k,j}$ case. The clique-based strategy outperforms our proposal in the constant $\delta^*_{k,j}$ case if $d_H$ is high enough, but only when the $n_k/p$ ratio is not too small. In the random $\delta^*_{k,j}$ case, our proposal clearly outperforms the clique-based strategy when $p=500$, and more generally as $d_H$ increases and/or the $n_k/p$ ratio decreases. These results suggest that the clique-based strategy might not be able to fully account for, nor benefit from, the true structure in a high-dimensional setting. They also suggest that our proposal is  better suited for this high-dimensional setting, as long as heterogeneity is not too high.


\section{Application on single-cell data}\label{sec:RealData}

We analyse data describing mielocytic leukemia cells undergoing differentiation to macrophage. Expression levels of 45 transcription factors are measured at $K=8$ distinct time points of this differentiation process  ($H0, H1, H6, H12, H24, H48, H72$ and $H96$). Each time point defines a stratum where data on $n_k = 120$ single cells are available. This data set is described in \cite{kouno2013temporal}. In this application, the main objective is to determine how associations among the 45 transcription factors vary over time.  \cite{kouno2013temporal} focus on marginal associations and use univariate analyses while graphical models, which describe conditional associations, might be better suited. Their inference can be reduced to the identification and description of the neighborhood of each covariate \citep{MeinBuhl2006}. Here, as a first step, we study how the neighborhood of  one particular transcription factor, EGR2, varies over time. Towards this end, we consider stratified linear regression models that relates EGR2 to the other $p=44$ factors on the $K=8$ strata. Expression levels of EGR2 are centered within each stratum, and no intercept is included in the models. Then, parameters of interest are vectors $\bbeta^*_1, \ldots, \bbeta^*_8$, where $\bbeta^*_k\in \R^{44}$ describes the association between EGR2 and the $p=44$ transcription factors at the $k$th time point. We compare estimates returned by our proposal and four competitors. The basic approach is considered with two distinct choices for the reference stratum: we set it to either $H0$ or $H96$ for each covariate.  We further consider the clique-based strategy. Finally, given the ordinal nature of the strata in this particular example, the variant based on chain graphs \citep{GertheissTutz} can be seen as the reference method. We include it as well, even if our main objective in this illustrative application is to compare the other four approaches, which do not account for this additional information. For each approach, regularization parameters are selected by 5-fold cross-validation. 

\begin{figure}[h]
\begin{center}
\includegraphics[scale=0.7]{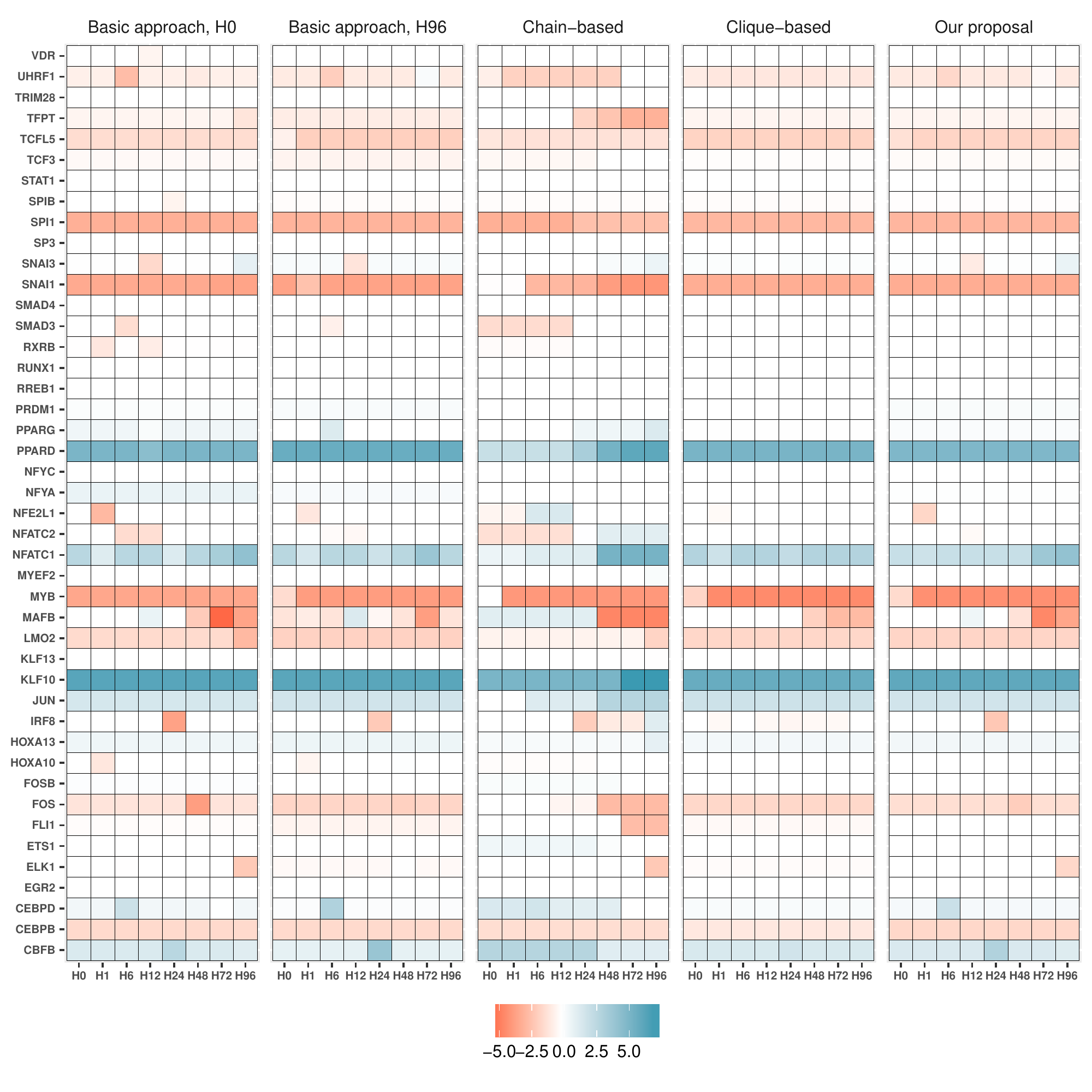}
\end{center}\caption{Estimation of the $K=8$ parameter vectors in the linear regression models describing the association between EGR2 and the $p=44$ other transcription factors, at times $H0, H1, H6, H12, H24, H48, H72$ et $H96$. Each column corresponds to the estimation obtained according to one of the five considered approaches : the basic approach with the reference stratum set to either $H0$ or $H96$ for every covariate, two versions of the generalized fused lasso estimates, one based on cliques and one based on chain graphs, and our proposal.}\label{Fig:EGR2}
\end{figure}

Results are presented in Figure \ref{Fig:EGR2}. For each method, the estimates correspond to a $44\times 8$ matrix of the form $(\hbbeta_1, \ldots, \hbbeta_K)$ with $K=8$ and $\hbbeta_k \in \R^{44}$ for any $k\in[K]$. Therefore, for any $j\in [p]$, the $j$th row of each matrix corresponds to estimates $(\hbeta_{1,j}, \ldots, \hbeta_{K,j})$ of the effect of the $j$th covariate over the $K$ time points, returned by the corresponding approach. This heat map representation allows an easy comparison of the pattern identified among the effects $(\hbeta_{1,j}, \ldots, \hbeta_{K,j})$ of each covariate across the $8$ strata. Our first objective is to illustrate the impact of the reference stratum when using the basic approach. Considering for instance the association between EGR2 and MYB, most approach identify the same pattern: the association is constant between $H1$ and $H96$, while it is lower, or even null, at $H0$.  However, setting the reference stratum to $H0$, the basic approach does not detect any heterogeneity.  These results are consistent with what is expected if the true pattern is the one identified by the other approaches: the basic approach used with $H0$ as the reference stratum is unlikely to identify the heterogeneity if it occurs at $H0$. Now consider the association with ELK1. Our proposal, the basic approach with $H0$ as the reference stratum and the strategy based on chain graphs all suggest the absence of association between EGR2 and ELK1 on the time interval $H0$ and $H72$, and a positive association at $H96$. If the reference stratum is set to $H96$, the basic approach suggests a quite different pattern, which is again expected if the true pattern is the one returned by the other approaches. These results confirm that the reference stratum is critical for the basic approach. They further suggest that our proposal is able to identify appropriate covariate-specific reference strata.

The comparison of the patterns returned by our approach and the two fused lasso strategies mostly highlights that the clique-based strategy identifies fewer heterogeneities and that strategy based on chain graphs returns smoother patterns. Again, these results were expected given the connectivity of the clique and chain graph, respectively. Prediction error was evaluated by double 5-fold cross-validation. Among the approaches that do not account for the ordering of the strata, the best prediction error is obtained with our proposal, while the worst is 1.8\% higher and is obtained with the clique-based strategy. The chain graph strategy leads to an improvement of 1.8\% compared to our approach. 

Two main conclusions can be drawn. When data come from several strata of the population and no information is available regarding which strata are likely to share similar effects, our proposal is a competitive approach. When additional information is available, as in this particular application where strata are naturally ordered, accounting for it can be beneficial.

\section{Discussion}\label{sec:Discussion}

After submitting a first version of this work, we became aware of concurrent work by \cite{GrossTibsh} where the authors introduce similar ideas. They apply it, in particular for the uplift problem in clinical research where the objective is to find sub-populations in a randomized trial for which an intervention is beneficial. In addition, there has been a recent line of works on penalized approaches aiming at identifying interactions, not necessarily between a categorical covariate and other predictors  \citep{LimHastie, Radchenko}. In this general context, strong hierarchy is often imposed: whenever an interaction between two variables is included in the model, the corresponding main effects are included too. However, this strong hierarchy is not desirable in our setting, where a coefficient can be nonzero in only one of the strata \citep{GrossTibsh}. Moreover, when applied in our setting, these approaches can be seen as versions of the basic approach (\ref{eq:RefLasso}) based on extensions of the $L_1$-norm penalty. In particular, a reference stratum has to be chosen as a first step, and the performance of these approaches would depend on this arbitrary choice. For instance, \cite{Radchenko} establish assumptions under which their approach is sparsistent. When applied with the reference strata $\ell$, their main condition on the design matrix is very similar to $(IC)_\ell$, which is generally stronger than $(IC)_{\ell^*}$ and $(IC)_0$. Then, these approaches could benefit from the ideas we developed in this article.

Our proposal is based on an overparametrization, which naturally raises the question of identifiability. We refer to \cite{GrossTibsh} for some discussion. We shall add that there is no identifiability issue under the conditions of Theorem \ref{theo:General}. If these conditions do not hold, and in particular if $(IC)_{\ell^*}$ is not fulfilled,  even the optimal version of the basic approach is not sparsistent, and the identifiability issue related to the overparametrization is secondary. 

Prediction bounds for our proposal can be derived under the conditions presented in this work. But weaker conditions might be sufficient, following recent work studying the lasso for correlated designs \citep{Dalalyan14}. Another extension might concern the derivation of valid p-values or confidence intervals for the nonzero parameters identified by our proposal. Given its connection with the lasso, this post-selection inference might be derived by extending recent strategies proposed for lasso estimates  \citep{LeeTaylor2016}. 

We also plan to extend our proposal to other regression models, which is straightforward for a variety of models given its connection with the lasso. In particular, our proposal could easily be extended to stratified Cox models used in survival analysis when competing risks arise \citep{rosner2013breast}, or to the conditional logistic models used in case-controls studies \citep{cLogitL1}. The extension of clique-based estimates to other models is generally more computationally burdensome, partly because there is no proximal operator for the fused penalty.

\section*{Acknowledgement}
We are grateful to Adeline Samson and Philippe Rigollet and the anonymous referees and associate editor for their fruitful comments on a preliminary version of this article. 

\section*{Supplementary material}
\label{SM}
Supplementary material available at {\em Biometrika} online includes two sections. Section 1 presents technical details: the proof of Lemma 1 and a generalized version of Lemma 1, the version of Theorem 1 in the balanced and orthogonal design along with its proof, and two corollaries describing the particular cases where $S_{\ell^*} = \emptyset$ and $T_{\ell^*}=\emptyset$. Section 2 presents additional results from our empirical study: accuracies for the recovery of other sets of interest in the settings described heres, and additional results obtained under an alternative settings which should favor the clique-based strategy.

\section*{Appendix: Supplementary Files}

\subsection{Technical details}

\subsubsection{Proof of Lemma 1}\label{proofLem1}
Lemma 1 is established for matrix $\bcX_0$; the proof for $\bcX_\ell$ follows from similar arguments and is omitted.

Fix $\tau_0>0$ and set $\tau_k=\tau_0 K^{-1/2}$ for all $k\in[K]$. 
Recall that $\btheta^{*}_0= ({\mu_{\ell^*}^{*}}^T,\tau_1{\bgamma_{0,1}^{*}}^T,\ldots,\tau_K{\bgamma_{0,K}^{*}}^T)^T$, with $\mu^*_{\ell^* j} = \beta^*_{\ell^*_j, j}$ for $j\in[p]$,  and $J_0 = \{j\in[(K+1)p]: \theta_{0 j}^{*}\neq 0\}.$ For the sake of brevity, the proof is only presented in the case where $S_{\ell^*}\neq \emptyset$ and $T_{\ell^*}\neq\emptyset$, where $S_{\ell^*}=\{j\in[p]: \mu^{*}_{\ell^* j} \neq 0\}$ and $T_{\ell^*}=\{(k,j): \beta^*_{k,j}\neq \mu^{*}_{\ell^* j} \}$. Setting $T^*_{k} = \{j\in[p]:(k,j)\in T_{\ell^*}\}$ and $\bSigma^{(k)}_{S_{\ell^*},T_{k}} = \bX_{S_{\ell^*}}^{(k)^T}\bX_{T_{k}}^{(k)}$ we have 
$$
\bcX_{0 J_0} = \left( 
\begin{array}{c c c c } 
\bX_{S_{\ell^*}}^{(1)} &\frac{\bX_{T^*_1}^{(1)}}{\tau_1}& \hdots &{0}  \\
\vdots & \vdots & \ddots & \vdots \\
\bX_{S_{\ell^*}}^{(K)} &{0}&\hdots  &  \frac{\bX_{T^*_k}^{(K)}}{\tau_K}  
\end{array}
\right),  
(\bcX_{0 J_0}^T\bcX_{0 J_0}) = \left( 
\begin{array}{c c c c c c }
n\bI_{|S_{\ell^*}|}&\frac{\bSigma^{(1)}_{S_{\ell^*},T_1^*}}{\tau_1}& \hdots & \hdots &\hdots & \frac{\bSigma^{(K)}_{S_{\ell^*},T_K^*}}{\tau_K} \\
\frac{\bSigma^{(1)^T}_{S_{\ell^*},T_1^*}}{\tau_1} & \frac{n\bI_{|T_1^*|}}{\tau_0^2}& {0}& \hdots & \hdots &{0} \\
\frac{\bSigma^{(2)^T}_{S_{\ell^*},T_2^*}}{\tau_2} &{0}&  \frac{n\bI_{|T_2^*|}}{\tau_0^2}& {0}&\hdots &{0}\\
\vdots &\vdots&\ddots & \ddots& \ddots & \vdots \\
\vdots &\vdots&\ddots & \ddots& \ddots &{0} \\
\frac{\bSigma^{(K)^T}_{S_{\ell^*},T_K^*}}{\tau_K}  & {0}&\hdots&\hdots&  {0}& \frac{n\bI_{|T_K^*|}}{\tau_0^2}
\end{array}
\right). 
$$
For any $s\in[|S_{\ell^*}|]$ and $t\in[|T^*_k|]$, denote the $s$th element of $S_{\ell^*}$ by $S_{\ell^* s}$  and the $t$th element of $T^*_k$ by $T^*_{k,t}$. For all $j\in S_{\ell^*}$, further denote by $N^*_j =n|K^*_j|/K$ the number of observations in strata contained in $K^*_j = \{k\in[K]:\beta^*_{k,j}= \beta^*_{\ell^*_j,j}\}$. Now introduce $\bC$, the matrix of size $|T_{\ell^*}|\times |S_{\ell^*}|$ made of $K$ blocks $\bC_k$. Each block is of size $|T^*_k|\times|S_{\ell^*}|$, and the element $(t,s)$ of $\bC_k$ is $(\bC_k)_{(t,s)} = -\tau_k/N^*_{S_{\ell^* s}}$ if $S_{\ell^* s}=T^*_{k,t}$ and 0 otherwise ($s\in[|S_{\ell^*}|]$ and $t\in[|T^*_k|]$). 

For $k,\ell\in[K]$ and $k\neq \ell$, further introduce ${B}_{k_1,k_2}$, the matrix of size $|T^*_{k_1}|\times |T^*_{k_2}|$ with element $(t_1,t_2)$ equal to $\tau_{k_1}\tau_{k_2}/N^*_{S_{\ell^* s}}$ if $T^*_{k_1,t_1}=T^*_{k_2,t_2}=S_{\ell^* s}$ for some $s\in[|S_{\ell^*}|]$, and 0 otherwise. For $k\in[K]$, denote by ${B}_{k,k}$ the diagonal matrix of size $|T^*_k|\times |T^*_k|$ with $t$th diagonal term equal to $\tau_k^2(N^*_{S_{\ell^* s}}+n_k)/(N_{\ell^* s}n_k)$ if $T^*_{k,t}=S_{\ell^* s}$ for some $s\in[|S^*|]$ and $\tau_k^2/n_k$ otherwise. Finally denote by ${D}$ the  diagonal matrix of size $|S_{\ell^*}|\times |S_{\ell^*}|$ with $j$th diagonal term equal to $1/N^*_j$, and by  $ {B}$ the matrix of size $|T_{\ell^*}|\times |T_{\ell^*}|$ made of $K^2$ blocks, with block $(k_1,k_2)$ equal to ${B}_{k_1,k_2}$. By standard algebra, we have
\begin{eqnarray*}
(\bcX_{0 J_0}^T\bcX_{0 J_0}) ^{-1} = \left( 
\begin{array}{l l}
{D}& \bC^T \\
\bC & {B}
\end{array}
\right).
\end{eqnarray*}
Now, for any $j\notin J_0$, the $j$th column $\bcX_{0 j}$ of $\bcX_0$ is of the form either (A) or (B):
\begin{itemize}
\item[(A)] $\bcX_{0 j} = ({0}^T, \ldots, {0}^T, X^{(k_0)^T}_{j_0}, {0}^T, \ldots, {0}^T)^T$ for some $k_0\in[K]$ and some $j_0\notin T^*_{k_0}$,
\item[(B)] $\bcX_{0 j} = (X^{(1)^T}_{j_0}, \ldots, X^{(K)^T}_{j_0})^T$ for some $j_0\notin S_{\ell^*}$.
\end{itemize}
If $\bcX_{0 j}$ is of form (A), then ${\bcX_{0 j}}^T\bcX_{0 J_0} = (\bar {d}^T, {0}_{|T_1^*|}^T, \ldots, {0}_{|T_K^*|}^T)\in\R^{|J_0|}$ with $\bar {d}\in\R^{|S_{\ell^*}|}$ and 
\[
\bar {d}_s = \left\{
\begin{array}{ll}
n_{k_0}/\tau_{k_0} & {\rm if}\ j_0 = S_{\ell^* s},\\
0 & {\rm otherwise}.
\end{array}\right.
\]
Therefore, if $\bcX_{0 j}$ is of form (A), we have 
\begin{equation*}
\|{\bcX_{0 j}}^T\bcX_{0 J_0} ({\bcX_{0 J_0}}^T\bcX_{0 J_0})^{-1} \|_1 \leq \max_{j\in S_{\ell^*}}\max_{k\in K^*_j} \frac{n_k}{\tau_k N^*_j}(1+ \sum_{\ell \notin K^*_j} \tau_\ell). 
\end{equation*}
If $\bcX_{0 j}$ is of form (B), then $\bcX_{0 j}^T\bcX_{0 J_0} = ({0}_{|S_{\ell^*}|}, {d}_1, \ldots, {d}_K)^T\in\R^{|J_0|}$, with  ${d}_k\in\R^{|T^*_k|}$ and 
\[
{d}_{k,t} = \left\{
\begin{array}{ll}
n_k/\tau_k & {\rm if}\ j_0 = T^*_{k,t},\\
0 & {\rm otherwise}.
\end{array}\right.
\]
In this case, we have 
\begin{equation*}
\|{\bcX_{0 j}}^T\bcX_{0 J_0} ({\bcX_{0 J_0}}^T\bcX_{0 J_0})^{-1} \|_1 \leq \max_{j\notin S_{\ell^*}}\sum_{\ell \notin K^*_j} \tau_\ell.
\end{equation*}

Moreover, under the setting considered in Lemma 1 we have $n_k=n/K$ for all $k\in[K]$ and $\tau_k=\tau_0 K^{-1/2}$. Therefore, $\max_{j\notin J_0}\|{\bcX_{0 j}}^T\bcX_{0 J_0} ({\bcX_{0 J_0}}^T\bcX_{0 J_0})^{-1} \|_1<1$ if and only if assumption $({sIC})_0$ holds, which completes the proof of Lemma 1 for matrix $\bcX_0$.

\subsubsection{Version of Theorem 1 in the setting of orthogonal designs and balanced strata}

Theorem \ref{Theo_SimpleCase} below is the version of Theorem 1 in the setting considered in Lemma 1, where $n_k = n/K$ and $(\bX^{(k)^T}\bX^{(k)})/n_k = \bI_{n_k}$ for all $k\in[K]$. We set ${\cal D}_0 = {\cal D}_{\ell^*, 0}$ and ${\cal D}_1 = {\cal D}_{\ell^*, 1}$; see Section 2.6 in the main text for the corresponding definitions.

\begin{theorem}\label{Theo_SimpleCase}
For all $k\in[K]$, assume that the noise variables $\varepsilon^{(k)}_i$, $i\in[n_k]$, are independent and identically distributed centered sub-Gaussian variables with parameter $\sigma>0$. Further assuming that  $({sIC})_0$ holds, define 
$$ \gamma = \min\left(1-{\cal D}_0\tau_0 K^{-1/2}, 1-\frac{K^{1/2} + {\cal D}_1\tau_0}{(K-{\cal D}_1)\tau_0}\right)$$
and $$C_{\min} = \displaystyle\min\left(1,\tau_0^{-2},\frac{1}{2}\left[ \Big(\tau_0^{-2} + 1\Big) - \left\{\Big(\tau_0^{-2} - 1\Big)^2 + \frac{4{\cal D}_1}{\tau_0^2 K}\right\}\right]\right). $$
For $\eta\in \{0,1\}$, we set
$$ \lambda^{(\eta)}_1 > \frac{2}{\gamma\min(1,\tau_0)}\left[\frac{2\sigma^2\log((K+\eta)p)}{n}\right]^{1/2}, \quad \lambda^{(\eta)}_{2,k} = \tau_k\lambda^{(\eta)}_1.$$
Finally introduce
$\beta_{\min}^{(\eta)}= \lambda^{(\eta)}_1[(|S_{\ell^*}| + |T_{\ell^*}|)^{1/2}C_{\min}^{-1} + 4\sigma C_{\min}^{-1/2}],$
and consider the following $\beta$-min conditions:
\begin{equation*}
({C}_{\beta^{(\eta)}_{\min}})(i):\ \forall j\in S_{\ell^*}, \ |\bar\beta^*_j| > \beta_{\min}^{(\eta)};\quad 
({C}_{\beta^{(\eta)}_{\min}})(ii):\ \displaystyle \forall j\in[p], \forall k\notin K^{*}_j, \ |\beta^*_{k,j} - \bar\beta^*_j| > \frac{K^{1/2}\beta_{\min}^{(\eta)}}{\tau_0}.
\end{equation*}
Then, $S_{\ell^*}$ and $T_{\ell^*}$ are both recovered 
\begin{itemize}
\vspace{-0.3cm}
\item with probability superior to $1- 4\exp(-c_1 n\lambda_1^{(0)^2})$, for some $c_1>0$, by the optimal version of the basic approach run with $\lambda_1 = \lambda_1^{(0)}$ and $\lambda_{2,k}=\lambda^{(0)}_{2,k}$ under $({C}_{\beta^{(0)}_{\min}})(i,ii)$ and we have $\|\widehat{\btheta}_{\ell^* J_{\ell^*}} - \btheta^{*}_{J_{\ell^*}}\|_\infty \leq \beta^{(0)}_{\min}$;
\item with probability superior to $1- 4\exp(-c_1 n\lambda_1^{(1)^2})$, for some $c_1>0$, by our approach run with $\lambda_1 = \lambda_1^{(1)}$ and $\lambda_{2,k}=\lambda^{(1)}_{2,k}$ under $({C}_{\beta^{(1)}_{\min}})(i,ii)$ and  we have $\|\widehat{\btheta}_{0 J_0} - \btheta^*_{0  J_0}\|_\infty \leq \beta^{(1)}_{\min}$.
\end{itemize}
\end{theorem}

Consider the asymptotic setting with $K$, and possibly $p$, tending to infinity as $n\rightarrow\infty$. Further assume that ${\cal D}_{\ell^*,0}={\cal D}_{\ell^*,1}={\cal D}_{\ell^*}$. If ${\cal D}_{\ell^*}\ll K^{1/2}$ or ${\cal D}_{\ell^*} = c K^{1/2}$ for some $0<c\leq 1/2$, then $\tau_0=1$ ensures perfect recovery for signals such that $\beta^{(\eta)}_{\min} = {\cal O}(n^{-1/2}[(|S_{\ell^*}| + |T_{\ell^*}|)\log((K+1)p)]^{1/2})$, which is optimal up to log-terms. If ${\cal D}_{\ell^*} = c K^{1/2}$ for some $c>1/2$, we get the same order of magnitude for $\beta^{(\eta)}_{\min}$, but with $\tau_0=(2c)^{-1}<1$. If $K^{1/2}\ll{\cal D}_{\ell^*}$, then the regime changes. For $K^{1/2}\ll{\cal D}_{\ell^*}\ll K$ the optimal choice is $\tau_0=K^{1/2}/(2{\cal D}_{\ell^*})$ which only ensures perfect recovery for $\beta^{(\eta)}_{\min} = {\cal O}((nK)^{-1/2}{\cal D}_{\ell^*}[(|S_{\ell^*}| + |T_{\ell^*}|)\log((K+1)p)]^{1/2})$. Finally, if ${\cal D}_{\ell^*}=cK$ for some $0<c<1/3$, then the result of Theorem S1 becomes almost meaningless: the optimal $\tau_0$ is ${\cal O}(K^{-1/2})$ which only ensures perfect recovery for $\beta^{(\eta)}_{\min} = {\cal O}(n^{-1/2}[K(|S_{\ell^*}| + |T_{\ell^*}|)\log((K+1)p)]^{1/2})$. 

\subsubsection{Proof of Theorem \ref{Theo_SimpleCase}}

In view of Lemma 1, and because ${\cal X}_{\ell^* J_{\ell^*}} =  \bcX_{0 J_0}$, we simply have to compute $\Lambda_{\min}((\bcX_{0 J_0}^T \bcX_{0 J_0}/n))$ in order to apply Theorem 1 of \cite{Wainwright2009} and establish Theorem \ref{Theo_SimpleCase}. 
To do so, we look for solutions of the characteristic polynomial of the matrix $(\bcX_{0 J_0}^T \bcX_{0 J_0}/n) $, $p(\lambda) = {\rm det}(  \bcX_{0 J_0}^T  \bcX_{0 J_0}/n - \lambda\bI_{|J_0|} )$. Using the block structure of matrix  $(  \bcX_{0 J_0}^T  \bcX_{0 J_0}/n - \lambda\bI_{|J_0|} )$, we get the following expression
\begin{equation*}
p(\lambda) = (1-\lambda)^{r_1} (\tau^{-2} - \lambda)^{r_2} \prod_{j\in S^*\cap\{\cup_k T_k^*\}} \left(\lambda^2 - \lambda(\tau^{-2}+1) +\tau^{-2}\Big(1 - \frac{n-N^*_j}{n}\Big)\right),
\end{equation*}
with $r_1 = |S_{\ell^*} \setminus \{\cup_k T^*_k\}|$ and $r_2 = |J_0| - |S_{\ell^*}| - |S_{\ell^*}\cap\{\cup_k T_k^*\}|$. It follows that 
$\Lambda_{\min}((  \bcX_{0 J_0}^T  \bcX_{0 J_0}/n))\geq \min\left(1, \tau^{-2},\frac{1}{2}\left[ (\tau^{-2} + 1) -\{(\tau^{-2} - 1)^2 + 4{\cal D}_1\tau^{-2}/ K\}^{1/2}\right]\right)$. Denote by $|\!\|M\|\!|_\infty$ the maximum row sum matrix norm of matrix $M$, and by $|\!\|M\|\!|_2$ its spectral norm. Because $|\!\|( \bcX_{0 J_0}^T \bcX_{0 J_0}/n)\|\!|_\infty \leq |J_0|^{1/2} |\!\|( \bcX_{0 J_0}^T \bcX_{0 J_0}/n)\|\!|_2 \leq (|S_{\ell^*}| + |T_{\ell^*}|)^{1/2}/C_{\min}$, Theorem \ref{Theo_SimpleCase} now follows from Theorem 1 of \cite{Wainwright2009}.

\subsubsection{Generalization of Lemma 1}

Here, we do not consider the orthogonal and balanced setting anymore and present general conditions ensuring that the irrepresentability conditions $(IC)_\ell$ and $(IC)_0$ are fulfilled by the design matrices $\bcX_\ell$ and $\bcX_0$ involved in the basic approach and our proposal, respectively. For all $k\in[K]$, we assume that $\tau_k = \tau_0(n_k/n)^{1/2}$ for some $\tau_0>0$ and that $n_k^{-1/2}\|X_j^{(k)}\|_2\leq 1$ for all $(j)\in[p]$. For $\bcX$ equal to either $\bcX_\ell$ or $\bcX_0$, this ensures that $n^{-1}\|{\cal X}_j\|_2\leq \max(1,\tau^{-1})$, for each column $\bcX_j$ of $\bcX$.

For any given vector of reference strata $\ell\in[K]^p$ and any $j\in[p]$, set $\bar K_{\ell,j} = \{k\in[K]: \beta^*_{k,j} = \beta^*_{\ell_j,j}\}$ and $K_{\ell,j} =\bar K_{\ell,j} \setminus\{\ell_j\}$. Further set, for any $k\in[K]$, $T_{\ell,k} = \{j\in [p]: \beta^*_{k,j}\neq \beta^*_{\ell_j, j}\}$,  $\bSigma_{\ell,k} = {\bX^{(k)T}_{T_{\ell,k}}}\bX^{(k)}_{T_{\ell,k}}$,  $
{{\Pi}_{\ell,k}} = \bX^{(k)}_{T_{\ell,k}} \bSigma_{\ell,k}^{-1}{\bX^{(k)T}_{T_{\ell,k}}}$, and $Z_{\ell,j}^{(k)} = (\bI_{n_k} - {\Pi}_{\ell,k}) X_j^{(k)}$. Define $\bomega_{\ell,j}^{(k)}= \bSigma_{\ell,k}^{-1}{\bX^{(k)^T}_{T_{\ell,k}}} X_j^{(k)}$ and $\bOmega_{\ell}^{(k)}= \bSigma_{\ell,k}^{-1}{\bX^{(k)^T}_{T_{\ell,k}}} \bX_{S_\ell}^{(k)}$. Introduce the quantities
$\widetilde \bSigma_{\ell} = \sum_{k\in[K]} {\bX_{S_\ell}^{(k)T}}(\bI_{n_k} - {\Pi}_{\ell,k})\bX_{S_\ell}^{(k)}$ and $\widetilde\bOmega_{\ell,j}^{(k)}= {\widetilde \bSigma_{\ell}}^{-1}{\bX_{S_\ell}^{(k)T}} Z_{\ell, j}^{(k)}$. Finally set 
\begin{eqnarray*}
c_{1}(\ell)&=&\max_{j\in S_\ell^{c}} \left\{\|\sum_{k\in [K]} \widetilde\bOmega_{\ell,j}^{(k)}\|_1 + \sum_{k\in [K]}\tau_k \|\sum_{l\in[K]} {\bOmega_{\ell}^{(l)}}\widetilde{\bOmega}_{\ell,j}^{(k)} \|_1 \right\} \\
c_{2}(\ell)&=&\max_{j\in[p]} \max_{k\in K_{\ell,j}} \left\{ \frac{\|\widetilde\bOmega_{\ell,j}^{(k)}\|_1}{\tau_k} + \sum_{l\neq k}\frac{\tau_l}{\tau_k}\|{\bOmega_{\ell}^{(l)}}\widetilde{\bOmega}_{\ell,j}^{(k)} \|_1 + \|{\bomega_{\ell,j}^{(k)}} + {\bOmega_{\ell}^{(k)}}  \widetilde{\bOmega}_{\ell,j}^{(k)}\|_1\right\}\\
\bar c_{2}(\ell)&=&\max_{j\in[p]} \max_{k\in \bar K_{\ell\!,j}}\!\! \left\{ \frac{\|\widetilde\bOmega_{\ell,j}^{(k)}\|_1}{\tau_k} + \sum_{l\neq k}\frac{\tau_l}{\tau_k}\|{\bOmega_{\ell}^{(l)}}\widetilde{\bOmega}_{\ell,j}^{(k)} \|_1 + \|{\bomega_{\ell,j}^{(k)}}   {\bOmega_{\ell}^{(k)}} + \widetilde{\bOmega}_{\ell,j}^{(k)}\|_1\right\}.
\end{eqnarray*}

\begin{lemma}\label{lem:IC.Ref}
Let $\ell\in[K]^p$ be a given vector of reference strata. Assume that $\Lambda_{\min}(\bSigma_{\ell,k})>0$ for $k\in[K]$ and $\Lambda_{\min}(\widetilde \bSigma_{\ell})>0$. Condition $(IC_\ell)$ holds  if and only if $c_1(\ell)<1$ and $c_2(\ell)<1$.

Assume that $\Lambda_{\min}(\bSigma_{\ell^*,k})>0$ for $k\in[K]$ and $\Lambda_{\min}(\widetilde \bSigma_{\ell^*})>0$. Condition $(IC)_0$ holds if and only if $c_1(\ell^*)<1$ and $\bar c_2(\ell^*)<1$.
\end{lemma}
The proof of Lemma \ref{lem:IC.Ref} follows from the same arguments as those presented in Section \ref{proofLem1} for the proof of Lemma 1, and is omitted.
Again, the conditions ensuring that $(IC)_{\ell^*}$ and $(IC)_0$ hold  are very similar. This shows that our proposal is able to mimic the optimal version of the basic approach even when the designs are not orthogonal or strata are not balanced.

\subsubsection{Other results}

Corollary \ref{Corol_Ident} and Corollary \ref{Corol_Indep}  consider the two special cases where $T_{\ell^*}=\emptyset$ and $S_{\ell^*}=\emptyset$, respectively.

\begin{corollary}\label{Corol_Ident}
Assume that the noise variables $(\varepsilon^{(k)}_i)_{i\in[n_k], k\in [K]}$ are independent and identically centered sub-Gaussian variables with parameter $\sigma>0$. Define $\mathbb{X} = (\bX^{(1)^T}, \ldots,\bX^{(K)^T})^T$, the $n\times p$ matrix with all the strata pooled together. Set $\tau_k=\tau_0(n_k/n)^{1/2}>0$ for all $k\in[K]$, for some $\tau_0>0$. For all $j\in[p]$, assume that there exists some $\beta^*_j\in\R$ such that $\beta^*_{k,j}=\beta^*_j$ for all $k\in[K]$ and set $S_{\ell^*}=\{j\in[p]: \beta^*_j\neq 0\}$. Further assume that $\Lambda_{\min}(\mathbb{X}^T_{S_{\ell^*}}\mathbb{X}_{S_{\ell^*}}/n)\geq C_{\min}$ for some $C_{\min} >0$, and that $n_k^{-1/2}\|X_j^{(k)}\|_2\leq 1$ for all $(k,j)\in[K]\times[p]$. Finally assume that the three following conditions hold:
\begin{eqnarray*}
&(\widetilde{\rm A})& \quad \sum_{k\in[K]}\tau_k > 1,\\
&(\widetilde{\rm C}.i.1)&\quad \tilde c_1:=\max_{j\notin S^{*}} \|(\mathbb{X}^T_{S_{\ell^*}}\mathbb{X}_{S_{\ell^*}})^{-1}\mathbb{X}^T_{S_{\ell^*}}\mathbb{X}_{j}\|_1< 1,\\
&(\widetilde{\rm C}.ii)&\quad \tilde c_2:=\max_{j\in[p]} \max_{k\in[K]} \tau_k^{-1}\| (\mathbb{X}^T_{S_{\ell^*}}\mathbb{X}_{S_{\ell^*}})^{-1} \bX^{(k)^T}_{S_{\ell^*}}X_j^{(k)}\|_1< 1.
\end{eqnarray*}
Now, set $\tilde\gamma = (1-\tilde c_1)\wedge(1-\tilde c_2)$ and 
$$ \lambda_1 = \frac{2}{(1\wedge \tau_0)\tilde\gamma}\left\{\frac{2\sigma^2\log((K+1)p)}{n}\right\}^{1/2},
\quad \beta_{\min}= \lambda_1\left(\frac{|S_{\ell^*}|^{1/2}}{C_{\min}} + 4\sigma C_{\min}^{-1/2} \right).$$
Then our proposal  run with parameters $\lambda_1$ and $\lambda_{2,k} = \tau_k \lambda_1$ identifies $S_{\ell^*}$ and $T_{\ell^*}=\emptyset$ with probability at least $1 - 4\exp(-c_1 n\lambda_1^2)$ for some $c_1>0$, as long as $\min_{j\in S^*} |\beta^*_{j}|> \beta_{\min}$.
\end{corollary}
Condition $(\widetilde{\rm C}.i.1)$ is exactly the irrepresentability condition on matrix $\mathbb{X}$, while conditions $(\widetilde{\rm C}.ii)$ and $(\widetilde{\rm A})$, which are very similar, both simply require that $\tau_0$ is high enough. Moreover, $\gamma = 1-\tilde c_1$ and $1\wedge \tau_0 = 1$ for $\tau$ high enough. Therefore, our proposal mimics the lasso run on $(\mathbb{X}, {\bcY})$ provided $\tau_0$ high enough and is optimal, up to log-terms, when the $\bbeta^*_{k}$'s are all equal. 

\begin{corollary}\label{Corol_Indep}
Assume that the noise variables $(\varepsilon^{(k)}_i)_{i\in[n_k], k\in [K]}$ are independent and identically centered sub-Gaussian variables with parameter $\sigma>0$. Set $\tau_k=\tau_0(n_k/n)^{1/2}>0$ for all $k\in[K]$, for some $\tau_0>0$. For all $j\in [p]$, set $K^{*}_j = \{k\in[K]: \beta^*_{k,j} =  0\}$ and for all $k\in[K]$, set $T^*_k = \{j\in[p]: \beta^*_{k,j}\neq 0\}$. Assume that $\min_{k}\{\Lambda_{\min}(\bX_{T^*_k}^{(k)^T}\bX_{T^*_k}^{(k)}/n_k)\}\geq C_{\min}$ for some $C_{\min}>0$, and that $n_k^{-1/2}\|X_j^{(k)}\|_2\leq 1$ for all $(k,j)\in[K]\times[p]$. Further assume that the three following conditions hold: 
\begin{eqnarray*}
&({\bar{\rm A}})& \quad \displaystyle \forall j\in[p], \sum_{\ell\notin K^{*}_j} \tau_\ell < 1 + \sum_{k\in K^{*}_j} \tau_k,\\
&({\bar{\rm C}}.i.1)&\quad {\bar{c}_1}:=\max_{j\in[p]} \max_{k\in[K]} \sum_{k\in[K]}\tau_k \| ({\bX^{(k)^T}_{T^*_k}}\bX^{(k)}_{T^*_k})^{-1}{\bX^{(k)^T}_{T^*_k}} X^{(k)}_j\|_1< 1,\\
&({\bar{\rm C}}.ii)&\quad {\bar{c}_2}:=\max_{k\in[K]}\max_{j\notin T^*_k}  \| ({\bX^{(k)^T}_{T^*_k}}\bX^{(k)}_{T^*_k})^{-1}{\bX^{(k)^T}_{T^*_k}} X^{(k)}_j\|_1< 1.
\end{eqnarray*}
Now, set ${\bar\gamma} = (1- {\bar c}_1)\wedge(1- {\bar c}_2)$,  and
$$ \lambda_1 = \frac{2}{(1\wedge \tau_0) {\bar\gamma}}\left\{\frac{2\sigma^2\log((K+1)p)}{n}\right\}^{1/2},
\quad \beta_{\min}= \lambda_1\left(\frac{\tau_0 |T_{\ell^*}|^{1/2}}{C_{\min}} +  4\sigma C_{\min}^{-1/2} \right).$$
Then our proposal run with parameters $\lambda_1$ and $\lambda_{2,k} = \tau_k \lambda_1$ recovers $S_{\ell^*}=\emptyset$ and $T_{\ell^*}=\{(k,j): \beta^*_{k,j}\neq 0\}$ with probability at least $1 - 4\exp(-c_1n\lambda_1^2)$, for some $c_1>0$, as long as $\min_{(k,j)\in T^*} |\beta^*_{k,j}|> \beta_{\min}(n/n_k)^{1/2}$.
\end{corollary}

Condition $({\bar{\rm C}}.ii)$ is exactly the union of the irrepresentability conditions for each matrix $\bX^{(k)}$, $k\in[K]$, while conditions $({\bar{\rm C}}.i.1)$ and $ {\bar{\rm A}}$ both simply require that $\tau_0$ is small enough. For $\tau_0$ small enough, our proposal then mimics the strategy consisting in performing $K$ lasso on the data $(\bX^{(k)},\by^{(k)})$, $k=1, \ldots, K$, independently, with a common  $\lambda_1$ value for each lasso. It is optimal, up to log-terms, when  $S_{\ell^*}=\emptyset$.

\subsection{Additional empirical results}

\subsubsection{Under the designs considered in the main text}

\begin{figure}[h]
\begin{center}
\includegraphics[scale=0.38]{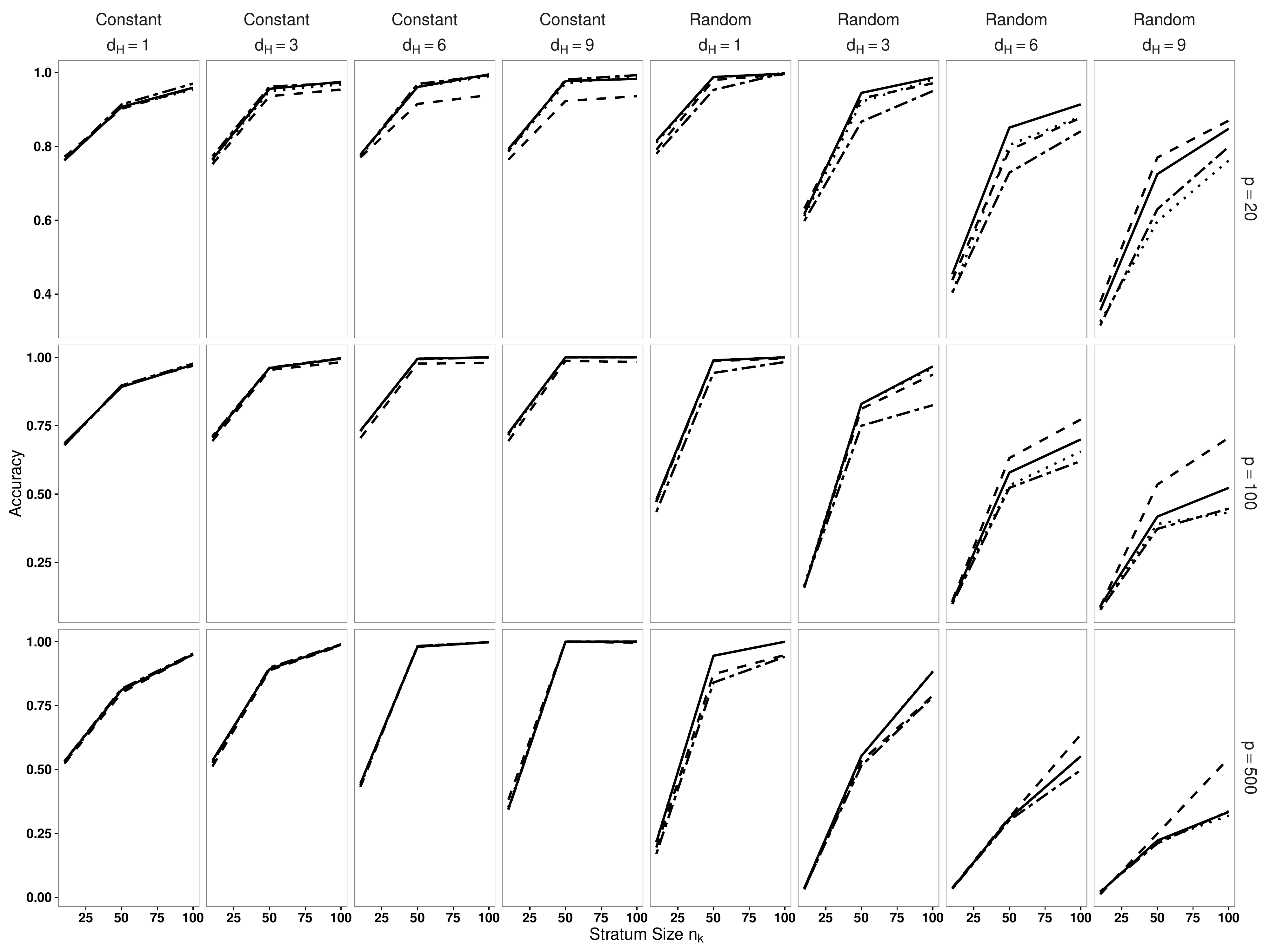}
\end{center}\caption{Accuracy regarding the recovery the set $S^{*}_{1, P_0}=\{j\in P_0: \beta^*_{1,j}\neq 0\}$ for the basic approach and of the set  $S^*_{P_0}=\{j\in P_0: \beta^*_{\ell^*_j,j}\neq 0\}$ for the three other approaches. {\em (Left)}: Constant $\delta^*_{k,j}$ case. {\em (Right)}: Random $\delta^*_{k,j}$ case. Results correspond to averages over 50 replicates in each configuration. Solid line: our proposal. Dotted line: optimal version of the basic approach. Dash-dot line: basic approach. Dashed line: clique-based approach.}\label{fig:AccS}
\end{figure}

Figure \ref{fig:AccS} presents additional results regarding the recovery of the set $S^*_{1, P_0}=\{j\in P_0: \beta^*_{1,j}\neq 0\}$ for the basic approach run with reference strata $\ell=(1, \ldots, 1)$ and the recovery of the set $S^*_{P_0}=\{j\in P_0: \beta^*_{\ell^*_j,j}\neq 0\}$ for the other approaches. Overall, our proposal and the optimal version of the basic approach perform similarly according to this criterion too. In the constant $\delta^*_{k,j}$ case, all methods perform similarly, and their performance does not depend on either $p$ or $d_H$. In the random $\delta^*_{k,j}$ case, the performance of each method decreases as $p$ and/or $d_H$ increases. This discrepancy with the results obtained in the constant case illustrates that it is harder for the optimal version of the basic approach, our proposal and the clique-based strategy to determine whether the overall effect $\beta^*_{\ell^*_j,j}$ of any covariate $j$ is null when the collection of values $(\beta^*_{1,j}, \ldots, \beta^*_{K,j})$ varies around zero; keep in mind that in the random case, we have either $\beta^*_{k,j}=1$ or $\beta^*_{k,j}=1 \pm \delta^*_{k,j}$ with $\delta^*_{k,j}\sim {\cal U}_{[K^{1/2}/2, 2K^{1/2}]}$ so that $\beta^*_{k,j}$ can be negative. Interestingly, it is harder for the basic approach to determine whether $\beta^*_{1,j}$ is null in this situation too. In addition, the basic approach is generally outperformed by the other approaches in this random $\delta^*_{k,j}$ case. The clique-based strategy performs well for $d_H\geq 6$, especially when  $n_k/p$ is not too small. But it is outperformed by our proposal, and the optimal version of the basic approach, for $d_H=3$ and $p=500$. 

\subsubsection{Under an alternative scenario}

\begin{figure}[h]
\begin{center}
\includegraphics[scale=0.36]{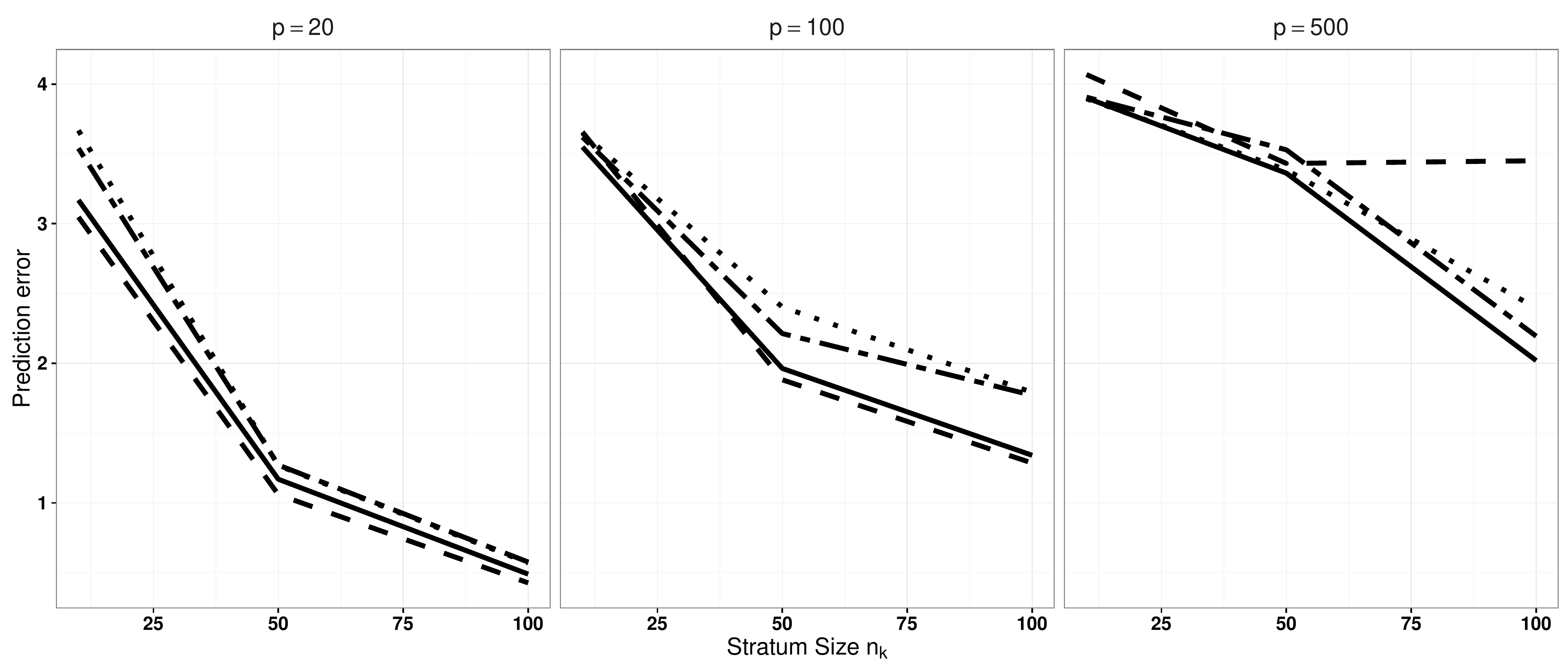}
\end{center}\caption{Prediction error (the lower, the better). Results correspond to averages over 50 replicates in each configuration. Solid line: our proposal. Dotted line: basic approach with $\ell=(20, \ldots, 20)$. Dash-dot line: basic approach with $\ell=(1, \ldots, 1)$. Dashed line: clique-based approach.}\label{fig:newSimL2Pred}
\end{figure}

\begin{figure}[h]
\begin{center}
\includegraphics[scale=0.5]{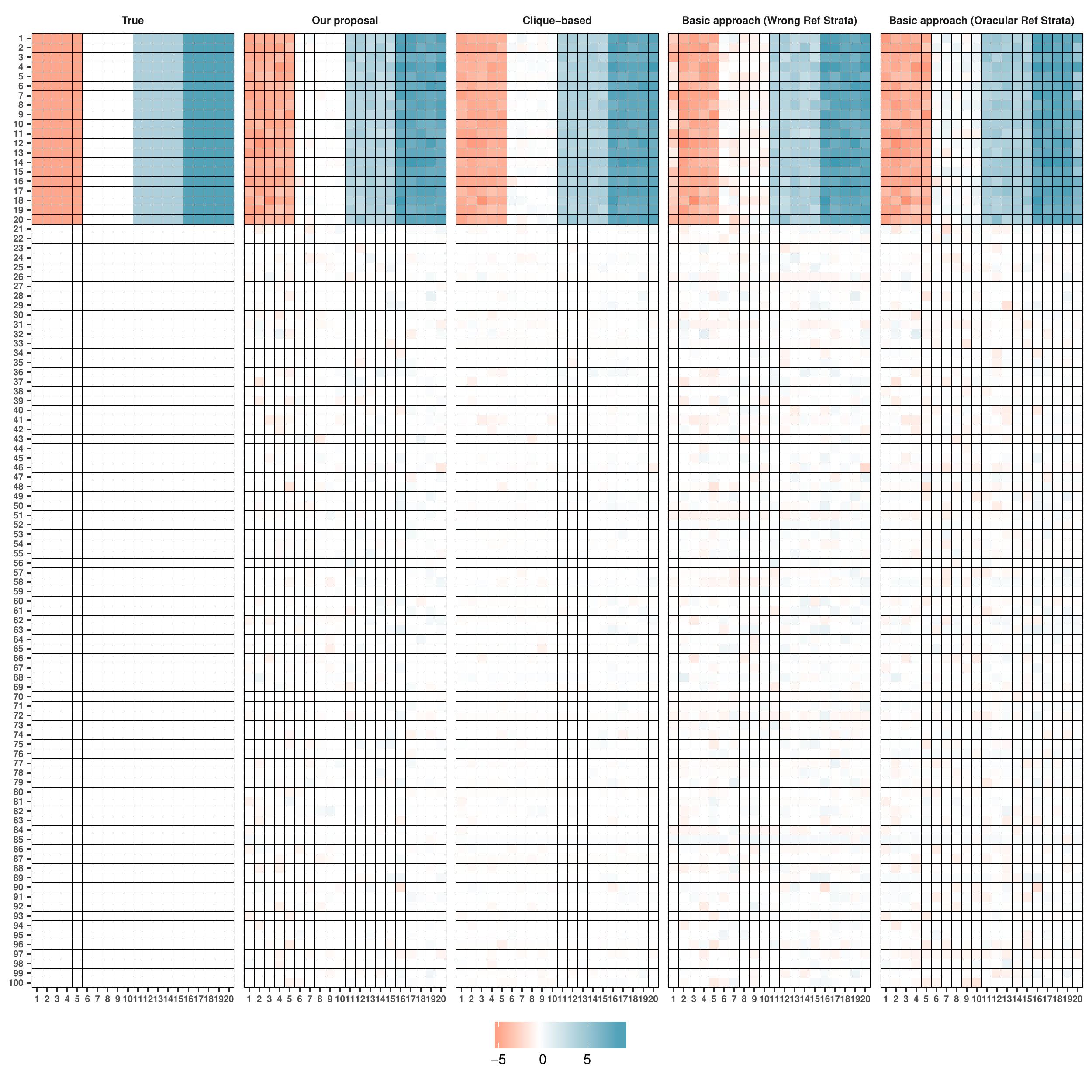}
\end{center}\caption{Estimation of the $K=20$ parameter vectors in one particular simulation with $n_k=100$ and $p=100$. The first column presents the true values. Each of the four remaining columns presents the estimates obtained according to one of the four considered approaches: our proposal, the clique-based approach and the basic approach with the reference stratum set to either $1$ or $20$ for every covariate.}\label{fig:HeatMapSim}
\end{figure}

Here, we present additional empirical results obtained under a scenario which should favor the clique-based strategy. We still consider the case where $K=20$ and $P_0\subset[p]$, with $|P_0|=20$ but, for each $j\in[P_0]$, we set $\beta_{1,j}=\cdots = \beta_{5,j} = -a$, $\beta_{6,j}=\cdots = \beta_{10,j} = 0$, $\beta_{11,j}=\cdots = \beta_{15,j} = a$ and $\beta_{16,j}=\cdots = \beta_{20,j} = 2a$, for some $a>0$. In other words, for each $j\in[P_0]$, the effects of the $j$th covariate across the 20 strata are made of 4 groups of distinct values, which should favor the clique-based strategy. 
Two values of $a$ were considered, $a=K^{1/2}$ and $a=K^{1/2}/3$. Because results were very similar, only those obtained for $a=K^{1/2}/3$ are presented here. Figure \ref{fig:newSimL2Pred} presents the predictive performance of each approach. 
We especially observe that our proposal performs nearly as well as the clique-based strategy in general, and outperforms it for $p=500$. It also slightly outperforms the two versions of the basic approach. 

Figure \ref{fig:HeatMapSim} presents the estimates returned by each approach for one particular simulation in the configuration $n_k=100$ and $p=100$. It can especially be seen that the clique-based strategy does not fuse coefficients sensibly better than the other approaches.

\end{document}